\documentclass[aop,MSNbibl,seceqn,dvips]{arximspdf}

\doi{10.1214/13-AOP877} 
\volume{42}
\issue{3}
\pubyear{2014}
\firstpage{1161}
\lastpage{1196}

\makeatletter

\newcommand{\rrvert}{\vert}
\newcommand{\llvert}{\vert}
\newcommand\Vol{{\operatorname{Vol}}}
\newcommand\R{{\mathbf{R}}}
\newcommand\C{{\mathbf{C}}}
\renewcommand\P{{\mathbf{P}}}
\newcommand\E{{\mathbf{E}}}
\newcommand\tr{{\operatorname{tr}}}
\newcommand\Z{{\mathbf{Z}}}
\newcommand\N{{\mathbf{N}}}
\newcommand\D{{\mathbf{D}}}
\newcommand\col{{\mathbf{c}}}
\newcommand\row{{\mathbf{r}}}
\newcommand\eps{\varepsilon}
\newcommand\cir{{\mathrm{cir}}}

\newcommand\Ba{\mathbf{a}}
\newcommand\Bf{\mathbf{f}}
\newcommand\Bu{\mathbf{u}}
\newcommand\Bv{\mathbf{v}}
\newcommand\Bx{\mathbf{x}}
\newcommand\By{\mathbf{y}}

\newcommand\BI{\mathbf{I}}

\newcommand\ep{\varepsilon}
\newcommand\1{\mathbf{1}}
\newcommand\Prob{\mathbf{P}}

\newcommand{\eqref}[1]{(\ref{#1})}

\newtheorem{theorem}[subsection]{Theorem}

\newproclaim{assumption}[subsection]{Assumption}

\newtheorem{lemma}[subsection]{Lemma}
\newtheorem{corollary}[subsection]{Corollary}

\newproclaim{example}[subsection]{Example}

\newtheorem{claim}[subsection]{Claim}

\newproclaim{definition}[subsection]{Definition}
\newproclaim{remark}[subsection]{Remark}
\newproclaim{remarks}[subsection]{Remarks}
\makeatother

\begin{document}
\begin{frontmatter}

\title{Random doubly stochastic matrices: The~circular~law}
\runtitle{Random doubly stochastic matrices: The circular law}

\begin{aug}
\author[A]{\fnms{Hoi H.} \snm{Nguyen}\corref{}\thanksref{m}\ead[label=e]{hoi.nguyen@yale.edu}}
\runauthor{H. H. Nguyen}
\affiliation{Ohio State University}
\thankstext{m}{Supported in part by NSF Grant DMS-12-00898.}
\address[A]{Department of Mathematics\\
Ohio State University\\
100 Math Tower\\
231 West 18th Avenue\\
Columbus, Ohio 43210\\
USA\\
\printead{e}} 
\end{aug}

\received{\smonth{4} \syear{2012}}
\revised{\smonth{8} \syear{2013}}

%
\begin{abstract}
Let $X$ be a matrix sampled uniformly from the set of doubly stochastic
matrices of size $n\times n$. We show that the empirical spectral
distribution of the normalized matrix $\sqrt{n}(X-\E X)$ converges
almost surely to the circular law. This confirms a conjecture of
Chatterjee, Diaconis and Sly.
\end{abstract}

%

\end{frontmatter}

\section{Introduction}\label{sectionintroduction}

Let $M$ be a matrix of size $n\times n$, and let $\lambda_1,\ldots
,\lambda_n$ be the eigenvalues of $M$. The empirical spectral
distribution (ESD) $\mu_{M}$ of $M$ is defined as
\[
\mu_{M}:=\frac{1}{n} \sum_{i\le n}
\delta_{\lambda_i}.
\]

We also define $\mu_{\cir}$ as the uniform distribution over the unit disk
\[
\mu_{\cir}(s,t):= \frac{1}{\pi} mes \bigl(|z|\le1; \Re(z)\le s, \Im
(z)\le t \bigr).
\]

Given a random $n \times n$ matrix $M$, an important problem in random
matrix theory is to study the limiting distribution of the empirical
spectral distribution as $n$ tends to infinity. We consider one of the
simplest random matrix ensembles, when the entries of $M$ are i.i.d.
copies of the random variable $\xi$.

When $\xi$ is a standard complex Gaussian random variable, $M$ can be
viewed as a random matrix drawn from the probability distribution $
\Prob(d M) =   \frac{1}{\pi^{n^2}} e^{- \tr(M M^\ast)} \,d M$
on the set of complex $n \times n$ matrices. This is known as the
\textit{complex Ginibre ensemble}. Following Ginibre \cite{Gi}, one may
compute the joint density of the eigenvalues of a random matrix $M$
drawn from the following complex Ginibre ensemble: $(\lambda_1, \ldots,
\lambda_n)$ has density
%
%
%
\begin{equation}
\label{eqginibre} p(z_1, \ldots, z_n):=
\frac{n!}{\pi^n \prod_{k=1}^n k! } \exp \Biggl( - \sum_{k=1}^n
|z_k|^2 \Biggr) \prod_{1 \leq i < j \leq n}
|z_i - z_j |^2
\end{equation}
on the set $|z_1| \leq\cdots\leq|z_n|$.

Mehta \cite{M,MB} used the joint density function \eqref{eqginibre}
to compute the limiting spectral measure of the complex Ginibre
ensemble. In particular, he showed that if $M$ is drawn from the
complex Ginibre ensemble, then $\mu_{({1}/{\sqrt{n}}) M}$ converges
to the circular law $\mu_{\cir}$. Edelman \cite{Ed-cir} verified the
same limiting distribution for the real Ginibre ensemble.

For the general case, there is no formula for the joint distribution of
the eigenvalues, and the problem appears much more difficult. The
universality phenomenon in random matrix theory asserts that the
spectral behavior of a random matrix does not depend on the
distribution of the atom variable $\xi$ in the limit $n \rightarrow
\infty$.

In the 1950s, Wigner \cite{W} proved a version of the universality
phenomenon for Hermitian random matrices. However, the random matrix
ensemble described above is not Hermitian; in fact, many of the
techniques used to deal with Hermitian random matrices do not apply to
non-Hermitian matrices.

An important result was obtained by Girko \cite{G1,G2} who related the
empirical spectral measure of non-Hermitian matrices to that of
Hermitian matrices. Consider the \textit{Stieltjes transform} $s_n$ of
$\mu_{({1}/{\sqrt{n})}M}$ given by
\[
s_n(z):= \frac{1}{n} \sum_{i=1}^n
\frac{1}{ ({1}/{\sqrt {n}})\lambda
_i- z} = \int_{\C} \frac{1}{ x + \sqrt{-1} y - z} \,d
\mu_{{1}/{\sqrt
{n}} M}(x,y).
\]

Since $s_n$ is analytic everywhere except at the poles, the real part
of $s_n$ determines the eigenvalues. We have
%
%
%
\begin{eqnarray}\label {eqrealstz}
\Re\bigl( s_n(z)\bigr) &=& \frac{1}{n} \sum
_{i=1}^n \frac{ ({1}/{\sqrt {n}}) \Re
(\lambda_i) - \Re(z) }{ \vert ( {1}/{\sqrt{n}}) \lambda_i - z
\vert ^2}
\nonumber
\\
&= &-\frac{1}{2n} \sum_{i=1}^n
\frac{\partial}{\partial\Re(z)} \log \biggl\llvert \frac{1}{\sqrt{n}}\lambda_i - z
\biggr\rrvert ^2
\\
&=& -\frac{1}{2n} \frac{\partial}{\partial\Re(z)} \log\det \biggl( \frac{1}{\sqrt{n}} M -
z \BI \biggr) \biggl( \frac{1}{\sqrt{n}} M - z\BI \biggr)^{\ast},
\nonumber
\end{eqnarray}
where $\BI$ denotes the identity matrix.

In other words, the task of studying the eigenvalues of the
non-Hermitian matrix $\frac{1}{\sqrt{n}}M$ reduces to studying the
eigenvalues of the Hermitian matrix $( \frac{1}{\sqrt{n}} M- z \BI) (
\frac{1}{\sqrt{n}} M - z\BI)^{\ast}$. The difficulty now is that the
$\log$ function has two poles, one at infinity and one at zero. The
largest singular value can easily be bounded by a polynomial in $n$.
The main difficulty is controlling the least singular value.

The first rigorous proof of the circular law for general distributions
was by Bai~\cite{B}. He proved the result under a number of moment and
smoothness assumptions on the atom variable $\xi$. Important results
were obtained more recently by Pan and Zhou \cite{PZ} and G\"otze and
Tikhomirov \cite{GT}. Using a strong lower bound on the least singular
value, Tao and Vu \cite{TVcir'} were able to prove the circular law
under the assumption that $\E|\xi|^{2+\eps} < \infty$, for some
$\eps>
0$. Recently, Tao and Vu (Appendix by Krishnapur) \cite{TVcir}
established the law assuming only that $\xi$ has finite variance.

\begin{theorem}[(\cite{TVcir})]\label{theoremcir}
Assume that the entries of $M$ are i.i.d. copies of a complex random
variable of mean zero and variance one, then the ESD of the matrix
$\frac{1}{\sqrt{n}}M$ converges almost surely to the circular measure
$\mu_\cir$.
\end{theorem}

In view of the universality phenomenon, it is important to study the
ESD of random matrices with nonindependent entries. Probably one of
the first results in this direction is due to Bordenave, Caputo and
Chafai \cite{BCC} who proved the following.

\begin{theorem}[(\cite{BCC}, Theorem~1.3)] \label{theoremBCC}
Let $X$ be a random matrix of size $n \times n$ whose entries are
i.i.d. copies of a nonnegative continuous random variable with finite
variance $\sigma^2$ and bounded density function. Then with probability
one the ESD of the normalized matrix $\sqrt{n}\bar{X}$, where $\bar
{X}=(\bar{x}_{ij})_{1\le i,j \le n}$ and ${\bar
{x}}_{ij}:=x_{ij}/(x_{i1}+\cdots+x_{in})$, converges weakly to the
circular measure $\mu_\cir$.
\end{theorem}

In particular, when $x_{11}$ follows the exponential law of mean one,
Theorem~\ref{theoremBCC} establishes the circular law for the
Dirichlet Markov ensemble (see also \cite{C}).

Related results with a linear assumption of independence include a
result of Tao, who among other things proves the circular law for
random zero-sum matrices.

\begin{theorem}[(\cite{Tao}, Theorem~1.13)]
Let $X$ be a random matrix of size $n \times n$ whose entries are
i.i.d. copies of a random variable of mean zero and variance one. Then
the ESD of the normalized matrix $\frac{1}{\sqrt{n}}\bar{X}$, where
$\bar{X}=(\bar{x}_{ij})_{1\le i,j \le n}$ and $\bar
{x}_{ij}:=x_{ij}-\frac{1}{n}(x_{i1}+\cdots+x_{in})$, converges almost
surely to the circular measure $\mu_\cir$.
\end{theorem}

With a slightly different assumption of dependence, Vu and the current
author showed in \cite{NgV-comb} the following.

\begin{theorem}[(\cite{NgV-comb}, Theorem~1.2)]\label{theoremNgV} Let
$0<\ep\le1$ be a positive constant. Let $M_n$ be a random $(-1,1)$
matrix of size $n\times n$ whose rows are independent vectors of
given\vadjust{\goodbreak}
row-sum $s$ with some $s$ satisfying $|s|\le(1-\ep)n$. Then the ESD of
the normalized matrix $\frac{1}{\sigma\sqrt{n}}M_n$, where $\sigma
^2=1-(\frac{s}{n})^2$, converges almost surely to the distribution
$\mu
_{\cir}$ as $n$ tends to $\infty$.
\end{theorem}

To some extent, the matrix model in Theorem~\ref{theoremNgV} is a
discrete version of the random Markov matrices considered in Theorem~\ref{theoremBCC}
 where the entries are now restricted to $\pm1/s$.
However, it is probably more suitable to compare this model with that
of random Bernoulli matrices. By Theorem~\ref{theoremcir}, the ESD of
a normalized random Bernoulli matrix obeys the circular law, and hence
Theorem~\ref{theoremNgV} serves as a local version of the law.

Although the entries of the matrices above are mildly correlated, the
rows are still independent. Because of this, we can still adapt the
existing approaches to bear with the problems. Our focus in this note
is on a matrix model whose rows and columns are not independent.

\begin{theorem}[(Circular law for random doubly stochastic
matrices)]\label{theoremMain}
Let $X$ be a matrix chosen uniformly from the set of doubly stochastic
matrices. Then the ESD of the normalized matrix $\sqrt{n}(X-\E X)$
converges almost surely to $\mu_\cir$.
\end{theorem}

Little is known about the properties of random doubly stochastic
matrices as it falls outside the scope of the usual techniques from
random matrix theory. However, there have been recent breakthrough by
Barvinok and Hartigan; see, for instance, \cite{BH1,BH2,BH3}. The
Birkhoff polytope $\mathcal{M}_n$, which is the set of doubly
stochastic matrices of size $n \times n$, is the basic object in
operation research because of its appearance as the feasible set for
the assignment problem. Doubly stochastic matrices also serve as a
natural model for priors in statistical analysis of Markov chains.
There is a close connection between the Birkhoff polytope and
$\operatorname{MS}(n,c)$, the set of matrices of size $n \times n$ with nonnegative
integer entries and all column sums and row sums equal $c$. These
matrices are called magic squares, which are well known in enumerative
combinatorics. We refer the reader to the work of Chatterjee, Diaconis
and Sly \cite{CDS} for further discussion.

There is a strong belief that random doubly stochastic matrices behave
like i.i.d. random matrices. This intuition has been verified in \cite
{CDS} in many ways. Among other things, it has been shown that the
normalized entry $nx_{11}$ converges in total variation to an
exponential random variable of mean one. More generally, the authors of
\cite{CDS} showed that the normalized projection $nX_k$, where $X_k$ is
the submatrix generated by the first $k$ rows and columns of $X$ with
$k=O(\frac{\sqrt{n}}{\log n})$, converges in total variation to the
matrix of independent exponential random variables.

Regarding the spectral distribution of $X$, it has been shown by
Chatterjee, Diaconis and Sly that the empirical distribution of the
singular values of $\sqrt{n}(X-\E X)$ obeys the quarter-circular law.
%
\begin{theorem}[(\cite{CDS}, Theorem~3)]\label{theoremsemi} Let $0\le
\sigma_1,\ldots,\sigma_n$ be the singular values of $\sqrt{n}(X-\E X)$,
where $X$ is a random doubly stochastic matrix. Then the empirical
spectral measure $\frac{1}{n}\sum_{i\le n}\delta_{\sigma_i}$ converges
in probability and in weak topology to the quarter-circle measure
$\frac
{1}{\pi}\sqrt{4-x^2}\mathbf{1}_{[0,2]} \,dx$.
\end{theorem}

The key ingredients in the proof of Theorem~\ref{theoremsemi} are a
sharp concentration result coupled with two transference principles
(Lemmas \ref{lemmarelation1} and \ref{lemmarelation2} below). These
principles help translate results from i.i.d. random matrices of
independent random exponential variables to random doubly stochastic matrices.

It has been conjectured in \cite{CDS} that the empirical spectral
distribution of $\sqrt{n}(X-\E X)$ obeys the circular law, which we
confirm now. For the rest of this section we sketch the general plan to
attack Theorem~\ref{theoremMain}.

Since the entries of $X$ are exchangeable, $\E X$ is the matrix $J_n$
of all $1/n$. The matrix $X-\E X$ has a zero eigenvalue, and we want to
single this outlier out due to several technical reasons. One way to do
this is passing to $\bar{X}$, a matrix of size $(n-1)\times(n-1)$
defined as
\[
\bar{X}:= \pmatrix{ x_{22}-x_{21} & \cdots&
x_{2n}-x_{21}\vspace*{2pt}
\cr
x_{32}-x_{31}
& \cdots& x_{3n}-x_{31}\vspace*{2pt}
\cr
\vdots& \vdots&
\vdots\vspace*{2pt}
\cr
x_{n2}-x_{n1} &\cdots&
x_{nn}-x_{n1} }.
\]

It is not hard to show that the spectra of $\sqrt{n}(X-\E X)$ is the
union of zero and the spectra of $\sqrt{n}\bar{X}$. Indeed, consider
the matrix $\lambda I_n - \sqrt{n}(X-\E X)$. By adding all other rows
to its first row, and then subtracting the first column from every
other column, we arrive at a matrix whose determinant is $\lambda\det
(\lambda I_{n-1}-\sqrt{n}\bar{X})$, thus confirming our observation.
Hence, it is enough to prove the circular law for~$\bar{X}$.

\begin{theorem}[(Main theorem)]\label{theoremmain}
Let $X$ be a matrix chosen uniformly from the set of doubly stochastic
matrices. Then the ESD of the matrix $\sqrt{n}\bar{X}$ converges almost
surely to $\mu_\cir$.
\end{theorem}

One way to prove our main result above is by showing that the Stieltjes
transform of $\mu_{\sqrt{n}\bar{X}}$ converges to that of the circular
measure. However, it is slightly more convenient to work with the
logarithmic potential. We will mainly rely on the following machinery
from \cite{TVcir}, Theorem~2.1.

\begin{lemma}\label{lemmareplacement}
Suppose that $M=(m_{ij})_{1\le i,j\le n}$ is a random matrix. Assume that:

\begin{itemize}
\item$\frac{1}{n}\|M\|_{HS}^2 =\frac{1}{n}\sum_{i,j}m_{ij}^2$ is
bounded almost surely;

\item for almost all complex numbers $z_0$, the logarithmic potential
$\frac{1}{n}\log|\det(M-z_0I_n)|$ converges almost surely to $f(z_0)=
\int_{\C} \log|w-z_0|\,d\mu_{\cir}(w)$.
\end{itemize}

Then $\mu_{M}$ converges almost surely to $\mu_{\cir}$.
\end{lemma}

We will break the main task into two parts, one showing the boundedness
and one proving the convergence.

\begin{theorem}\label{theorembounded}
Let $X$ be a matrix chosen uniformly from the set of doubly stochastic
matrices. Then there exists a constant $C$ such that the following holds:
\[
\P\biggl(\sum_{2\le i,j\le n} (x_{ij}-x_{i1})^2
\ge C\biggr) = O(\exp\bigl(-\Theta (\sqrt{n})\bigr).
\]
\end{theorem}

The proof of Theorem~\ref{theorembounded} will be presented at the end
of Section~\ref{sectionproperties}. The heart of our paper is to
establish the convergence of $\frac{1}{n}\log|\det(\sqrt{n}\bar
{X}-z_0I_{n-1})|$.

\begin{theorem}\label{theoremconvergence}
For almost all complex numbers $z_0$, $\frac{1}{n}\log|\det(\sqrt {n}\bar{X}-z_0I_{n-1})|$ converges almost surely to $f(z_0)$.
\end{theorem}

The main difficulty in establishing Theorem~\ref{theoremconvergence}
is that the entries in each row and each column of $\bar{X}$ are not at
all independent. To the best of our knowledge, the convergence for such
a model has not been studied before in the literature. We will present
its proof in Section~\ref{sectionproof}.

\textit{Notation}. Here and later, asymptotic notation such as
$O,\Omega,\Theta$ and so forth, are used under the assumption that
$n\rightarrow\infty$. A notation such as $O_C(\cdot)$ emphasizes that the
hidden constant in $O$ depends on $C$.

For a matrix $M$, we use the notation $\row_i(M)$ and $\col_j(M)$ to
denote its $i$th row and $j$th column, respectively. For an event $A$,
we use the subscript $\P_{\Bx}(A)$ to emphasize that the probability
under consideration is taking according to the random vector~$\Bx$.

For a real or complex vector $\Bv=(v_1,\ldots,v_n)$, we will use the
shorthand $\|\Bv\|$ for its $L_2$-norm $ (\sum_i|v_i|^2 )^{1/2}$.

\section{Some properties of random doubly stochastic matrices}\label
{sectionproperties}

We will gather here some basic properties of random doubly stochastic
matrices. The reader is invited to consult \cite{CDS} for further
insights and applications.

\subsection{Relation to random i.i.d. matrix of exponentials} Let
$\mathcal{M}_n$ be the Birkhoff polytope generated by the permutation
matrices. Let $\Phi$ be the projection from $\R^{n^2}$ to $\R
^{(n-1)^2}$ by mapping $(x_{ij})_{1\le i,j \le n}$ to $(x_{ij})_{2\le
i,j \le n}$.

Let $\Gamma\dvtx  \R^{(n-1)^2}\to\R^{n^2}$ denote the following function:
\[
\Gamma(X)=\Gamma(X)_{ij}:= \cases{ \displaystyle x_{ij}, & \quad$2\le i,j\le
n;$\vspace*{2pt}
\cr
\displaystyle 1-\sum_{k=2}^{n}x_{ik},
& \quad$2\le i\le n, j=1;$\vspace*{2pt}
\cr
\displaystyle 1-\sum_{k=2}^{n}x_{kj},
&\quad$ 2\le j\le n, i=1;$\vspace*{2pt}
\cr
\displaystyle 1-\sum_{l=2}^{n}
\Biggl(1-\sum_{k=2}^{n}x_{kl}
\Biggr), &\quad  $i=j=1.$}
\]

Thus $\Gamma$ extends a matrix $X$ of size $(n-1) \times(n-1)$ to a
doubly stochastic matrix of size $n\times n$ whose bottom right corner
is $X$. With the above notation, the doubly stochastic matrices
correspond to $(n-1)\times(n-1)$-matrices of the set
\[
S_n:= \bigl\{ X=(x_{ij})_{2\le i,j\le n}
\in[0,1]^{(n-1)^2}\dvtx 0\le \Gamma (X)_{ij}\le1 \bigr\}.
\]

The distribution of $X$ as a random doubly stochastic matrix is then
given by the uniform distribution on $S_n$. We next introduce an
asymptotic formula by Canfield and Mckay \cite{CM} for the volume of $S_n$,
%
%
\begin{equation}
\label{eqnCM} \Vol(S_n)=\frac{1}{n^{n-1}}\frac{1}{(2\pi)^{n-1/2}n^{(n-1)^2}}\exp
\biggl(\frac
{1}{3}+n^2+o(1)\biggr).
\end{equation}

This formula plays a crucial role in the transference principles to be
introduced next.

Define
\[
D_n:= \biggl\{ Y=(y_{ij})_{1\le i,j\le n}\dvtx \Phi\biggl(
\frac{1}{n}Y\biggr) \in S_n, \min \biggl\{\frac{1}{n}y_{ij}
-\Gamma\biggl(\Phi\biggl(\frac{1}{n}Y\biggr)\biggr)_{ij} \biggr\}
\ge0 \biggr\},
\]
where $\Phi\dvtx  \R^{n^2} \to\R^{(n-1)^2}$ is the projection
$X=(x_{ij})_{1\le i,j\le n} \mapsto(x_{ij})_{2\le i,j\le n}$.

Let $Y=(y_{ij})_{1\le i,j\le n}$ be a random matrix where $y_{ij}$ are
i.i.d. copies of a random exponential variable with mean one. As an
application of \eqref{eqnCM}, it is not hard to deduce the following
transference principle between random doubly stochastic matrices $X$
and random i.i.d. matrices $Y$.

\begin{lemma}[(\cite{CDS}, Lemma~2.1)]\label{lemmarelation1} Conditioning
on $Y\in D_n$, we have that $(\frac{1}{n}y_{ij})_{2\le i,j\le n}$ is
uniform on $S_n$. Furthermore, for large $n$ we have
\[
\P(Y\in D_n)\ge n^{-4n}.
\]
\end{lemma}

Lemma~\ref{lemmarelation1} is useful when we want to pass an
extremely rare event from the model $\frac{1}{n}Y$ to the model $X$. In
applications (in particular when working with concentration results),
it is more useful to work with matrices of bounded entries. With this
goal in mind we define
\[
\tilde{S}_n:= \biggl\{\tilde{X}=(\tilde{x}_{ij})_{2\le i,j \le n}
\in [0,1]^{(n-1)^2}, 0\le\Gamma(\tilde{X})_{ij} \le
\frac{10\log n}{n} \biggr\}
\]
and
\begin{eqnarray*}
&&\tilde{D}_n: = \biggl\{\tilde{Y}=(\tilde{y}_{ij})_{1\le i,j\le n}
\in [0,10\log n]^{n^2}, \frac{1}{n}\tilde{Y} \in
\tilde{S}_n, \\
&&\hspace*{78pt} 0\le\frac
{1}{n}\tilde{y}_{ij}-\Gamma
\biggl(\Phi\biggl(\frac{1}{n}\tilde{Y}\biggr)\biggr)_{ij}\le
n^{-4} \biggr\}.
\end{eqnarray*}

Observe that $\tilde{S}_n$ corresponds to doubly stochastic matrices
$\tilde{X}$ with entries bounded by $10\log n/n$.

Let $\tilde{Y}=(\tilde{y}_{ij})_{1\le i,j\le n}$ where $\tilde{y}_{ij}$
are i.i.d. copies of a truncated exponential $\tilde{y}$ with the
following density function:
%
%
\begin{equation}
\label{eqntilde{y}} \rho_{\tilde{y}}(x)= \cases{ \exp(-x)/
\bigl(1-n^{-10}\bigr), & \quad$\mbox{if } x\in[0,10\log n],$\vspace*{2pt}
\cr
0,
&\quad $\mbox{otherwise.}$}
\end{equation}

It is clear that $\E(\tilde{y}^2)=\Theta(1)$ and $\E(\tilde
{y}^4)=\Theta
(1)$. We now introduce another transference principle which is an
analogue of Lemma~\ref{lemmarelation1}.

\begin{lemma}[(\cite{CDS}, Lemma~4.1)]\label{lemmarelation2} Conditioning
on $\tilde{Y}\in\tilde{D}_n$, we have that $(\frac{1}{n}\tilde
{y}_{ij})_{2\le i,j\le n}$ is uniform on $\tilde{S}_n$. Furthermore,
for large $n$ we have
\[
\P(\tilde{Y}\in\tilde{D}_n)\ge n^{-10n}.
\]
\end{lemma}

Notice that in the corresponding definition of $\tilde{D}_n$ in \cite{CDS}, Section~4, the bound $10\log n$ was replaced by $6\log n$, but
one can easily check that this modification does not affect the
validity of Lemma~\ref{lemmarelation2}.

\subsection{Relation to random stochastic matrices} Let $\mathcal
{R}=\mathcal{R}_{r,n}$ denote the $r(n-1)$-dimensional polytope of
nonnegative matrices of size $r \times n$ whose rows sum to 1. Let $\mu
_r$ denote the uniform probability measure on $\mathcal{R}$, and let
$\nu_r$ denote the measure on $\mathcal{R}$ induced by the first $r$
rows of a random doubly stochastic matrix~$X$. As another application
of \eqref{eqnCM} (to be more precise, we need a more general form for
volume of polytopes generated by rectangular matrices of constant row
and column sums), one can show that these two measures are comparable
as long as $r$ is small.

\begin{lemma}[(\cite{CDS}, Lemma~3.3)]\label{lemmaderivative}
For a fixed integer $r\ge1$ and $n>r$ the Radon--Nikodym derivative of
the measures $\mu_r$ and $\nu_r$ satisfies
\[
\frac{d\nu_r}{d\mu_r} \le\bigl(1+o(1)\bigr)\exp(r/2)
\]
as $n\rightarrow\infty$.
\end{lemma}

It then follows that, in terms of order, there is not much difference
between the models $X$ and $\tilde{X}$.

\begin{theorem}\label{theorembound}
Assume that $B>4$ is a constant, then
\[
\P_X\bigl(n^{-B}\le nx_{11}\le B \log n\bigr)
\ge1- O\bigl(n^{-B/2}\bigr).
\]
\end{theorem}

In particular, since the entries of $X$ are exchangeable, Theorem~\ref
{theorembound} yields the following.

\begin{corollary}\label{corbound}
Assume that $X$ is a random doubly stochastic matrix, then
\[
\P(X\in\tilde{S}_n)= \P\bigl(|x_{ij}|\le10 \log n/n \mbox{ for
all } 1\le i,j\le n\bigr) \ge1-O\bigl(n^{-3}\bigr).
\]
\end{corollary}

\begin{pf*}{Proof of Theorem~\ref{theorembound}} It follows from Lemma~\ref
{lemmaderivative} (for $r=1$) that
\[
\P\bigl(n^{-B}\le nx_{11}\le B\log n\bigr) \le \bigl(1+o(1)
\bigr)\exp(1/2) \P \bigl(n^{-B}\le n x_1 \le B\log n\bigr),
\]
where $x_{1}$ has distribution $B(1,n-1)$.

The claim then follows because
\begin{eqnarray*}
&&\P\bigl(n^{-B}\le n x_1 \le B\log n\bigr)\\
&&\qquad=(n-1)\int
_{n^{-B}}^{B\log n} (1-x)^{n-2}\,dx
\\
&&\qquad=1-(n-1) \biggl(\int_{0}^{n^{-B}}
(1-x)^{n-2}\,dx + \int_{B\log n}^{1}
(1-x)^{n-2}\,dx \biggr)
\\
&&\qquad \ge1-O\bigl(n^{-B/2}\bigr).
\end{eqnarray*}
\upqed\end{pf*}

We end this section by giving a proof for the boundedness of Lemma~\ref
{lemmareplacement}.

\subsection{\texorpdfstring{A proof for Theorem \protect\ref{theorembounded}}
{A proof for Theorem 1.9}} We first
focus on the random vector $\Bx=(x_1,\ldots,x_n)$ chosen uniformly from
the simplex $S= \{\Bx=(x_1,\ldots,x_n), 0\le x_i\le1, \sum_i x_i
=1 \}$. Because each $x_i$ has distribution $B(1,n-1)$, we have
%
%
\begin{equation}
\label{eqnsimplex1} \E_{\Bx} \|\Bx\|^2 =
\frac{2}{n+1}.
\end{equation}

Also, it can be shown that (e.g., from \cite{MM}, equation (19))
%
%
\begin{equation}
\label{eqnsimplex2} \E_{\Bx} x_1x_2 =
\frac{1}{n(n+1)}.
\end{equation}

It thus follows from \eqref{eqnsimplex1} that $\|\Bx\|=O(1/\sqrt{n})$
with high probability. It turns out that this probability is extremely
close to one.

\begin{lemma}\label{lemmaisotropic}
Assume that $\Bx$ is sampled uniformly from $S$ and assume that $\ep>0$
is a sufficiently small constant. Then there exists a positive constant
$C>0$ such that
\[
\P\bigl(\|\Bx\| \ge C/\sqrt{n}\bigr) \le\exp(-\ep\sqrt{n}).
\]
\end{lemma}

We assume Lemma~\ref{lemmaisotropic} for the moment.

\begin{pf*}{Proof of Theorem~\ref{theorembounded}} First, it follows from
Lemma~\ref{lemmaderivative} (for $r=1$) that
\begin{eqnarray*}
&&\P\bigl(x_{21}^2+\cdots+x_{n1}^2\ge
C/n\bigr)\\
 &&\qquad\le \bigl(1+o(1) \bigr)\exp(1/2) \P \bigl(x_2^2+
\cdots+x_n^2 \ge C/n\bigr)
\\
&&\qquad=O(1)\P\bigl(x_1^2+x_2^2+
\cdots+x_n^2 \ge C/n\bigr),
\end{eqnarray*}
where $(x_1,x_2,\ldots,x_n)$ are sampled uniformly from the simplex $S$.
But Lemma~\ref{lemmaisotropic} indicates that the RHS is bounded by
$\exp(-\ep\sqrt{n})$. Thus
%
%
\begin{equation}
\label{eqnbound2} \P\bigl(x_{21}^2+\cdots+x_{n1}^2
\ge C/n\bigr)=O\bigl(\exp(-\ep\sqrt{n})\bigr).
\end{equation}

And so, as $x_{ij}$ are exchangeable, for any $j$ we also have
%
%
\begin{equation}
\label{eqnboundi} \P\bigl(x_{2j}^2+\cdots+x_{nj}^2
\ge C/n\bigr)=O\bigl(\exp(-\ep\sqrt{n})\bigr).
\end{equation}

The claim of Theorem~\ref{theorembounded} then follows because $\sum_{2\le i,j\le n}(x_{ij}-x_{i1})^2 \ge C$ would
imply $\sum_{i=2}^nx_{ij}^2\ge C/4n$ for some $j$.
\end{pf*}

It remains to prove Lemma~\ref{lemmaisotropic}. We show that it is a
direct consequence of the following geometric result.

\begin{theorem}[(\cite{P}, Theorem~1.1)]\label{theoremiso}
There exists an absolute constant $c>0$ such that if $K$ is an
isotropic convex body in $\R^n$, then
\[
\P\bigl(\Bx\in K, \|\Bx\| \ge c\sqrt{n}L_K t\bigr)\le\exp(-\sqrt{n}t)
\]
for every $t\ge1$, where $L_K$ is the isotropic constant of $K$.
\end{theorem}

Observe that, by the triangle inequality, for Lemma~\ref
{lemmaisotropic} it is enough to give a similar probability bound for
the event $\|\Bx-(1/n,\ldots,1/n)\|\ge C/\sqrt{n}$.

We first shift $S$ to the hyperplane $H:=\{\Bx'=(x_1',\ldots,x_n'),
x_1'+\cdots+x_n'=0\}$ by the translation $\Bx=(x_1,\ldots,x_n) \mapsto
(x_1-1/n,\ldots, x_n-1/n)$. We then scale the obtained body by a factor
$\alpha=\Theta(n)$ to obtain a regular simplex $S'$ of volume one.
Elementary computations show that this is an isotropic body of bounded
isotropic constant. Indeed, if $\Bx'=(x_1',\ldots,x_n')$ is sampled
uniformly from $S'$ and if $\bolds{\Theta}=(\theta_1,\ldots,\theta_n)$
is any unit vector in $H$, then by \eqref{eqnsimplex1} and \eqref
{eqnsimplex2},
\begin{eqnarray*}
\E_{\Bx'\in S'}\bigl\langle\Bx',\bolds{\Theta} \bigr
\rangle^2 &=& \E_{\Bx
'\in
S'} \biggl(\sum
_i \theta_i x_i'
\biggr)^2=\E_{\Bx\in S} \sum_i
\alpha^2 \biggl(\sum_i\theta
_i\biggl(x_i-\frac{1}{n}\biggr)
\biggr)^2
\\
&= &\alpha^2\sum_{i}
\theta_i^2 \biggl(x_i-\frac{1}{n}
\biggr)^2 + 2\alpha^2 \sum_{i
\neq j}
\theta_i\theta_j\biggl(x_i-
\frac{1}{n}\biggr) \biggl(x_j-\frac{1}{n}\biggr)
\\
&=& \alpha^2 \biggl(\frac{2}{n(n+1)}-\frac{1}{n^2}\biggr)\sum
_i \theta_i^2 + 2\alpha
^2 \biggl(\frac{1}{n(n+1)}-\frac{1}{n^2}\biggr)
\theta_i\theta_j
\\
&=&\alpha^2\biggl(\frac{1}{n(n+1)}\biggr) \sum
_i \theta_i^2 +
\alpha^2 \biggl(\frac
{1}{n(n+1)}-\frac{1}{n^2}\biggr) \biggl(
\sum_i \theta_i\biggr)^2
\\
&=&\frac{\alpha^2}{n(n+1)}.
\end{eqnarray*}

Thus the isotropic constant of $S'$ is of constant order. Theorem~\ref
{theoremiso} applied to $\Bx'$ yields the following for a sufficiently
large constant $C$:
\[
\P\bigl(\Bx'\in S', \bigl\|\Bx'\bigr\|\ge C
\sqrt{n}\bigr)\le\exp(-\ep\sqrt{n}).
\]
Lemma~\ref{lemmaisotropic} then follows because $\alpha\|\Bx
-(1/n,\ldots,1/n)\| = \|\Bx'\|$.

\begin{remark} The proof above heavily relies on the isotropic property
of the simplex $S$. It is perhaps more natural to relate $\Bx
=(x_1,\ldots
,x_n)$ to $(y_1/(y_1+\cdots+y_n),\ldots, y_n/(y_1+\cdots+y_n))$, where
$y_i$ are i.i.d. copies of a random exponential random variable of mean
one.\setcounter{footnote}{1}\footnote{The author is grateful to the anonymous referee for this
suggestion.} The probability $\P(\|\Bx\|\ge C/\sqrt{n})$ is then
bounded by the sum $\P(y_1+\cdots+y_n \le n/\sqrt{C}) + \P
(y_1^2+\cdots
+y_n^2 \ge Cn)$. As now we only need to work with sum of i.i.d. random
variables, by choosing $C$ sufficiently large, it is not hard to show
that both $\P(y_1+\cdots+y_n \le n/\sqrt{C})$ and $\P(y_1^2+\cdots+y_n^2
\ge Cn)$ are bounded from above by $\exp(-\Theta(\sqrt{n}))$.
\end{remark}

\section{\texorpdfstring{The singularity of $\bar{X}$}{The singularity of X}}
\label{sectionsingularity}
In order to justify Theorem~\ref{theoremconvergence}, one of the key
steps is to bound the singularity probability of the matrix $\sqrt {n}\bar{X}-z_0I_{n-1}$. This problem is of interest on its own.

We will show the following general result regarding the least singular
value~$\sigma_{n-1}$.

\begin{theorem}\label{theoremsingular}
Let $F=(f_{ij})_{2\le i,j\le n}$ be a deterministic matrix where
$|f_{ij}|\le n^\gamma$ with some positive constant $\gamma$. Let $X$ be
an $n\times n$ matrix chosen uniformly from the set of doubly
stochastic matrices. Then for any positive constant $B$ there exists a
positive constant $A$ such that
\[
\P\bigl(\sigma_{n-1}(\bar{X}+F)\le n^{-A}\bigr) \le
n^{-B}.
\]
\end{theorem}

Combine with Theorem~\ref{corbound} we obtain the following important
corollary which we reserve for later applications.

\begin{corollary}\label{corsingulartilde{X}} Let $F=(f_{ij})_{2\le
i,j\le n}$ be a deterministic matrix where $|f_{ij}|\le n^\gamma$ with
some positive constant $\gamma$. Let $\tilde{X}=(x_{ij})$ be a random
doubly stochastic matrix where $x_{ij}\le10\log n/n$ for all $1\le
i,j\le n$. Then there exists a positive constant $A$ such that
\[
\P\bigl(\sigma_{n-1}(\bar{\tilde{X}}+F)\le n^{-A}\bigr)= O
\bigl(n^{-3}\bigr).
\]
\end{corollary}

Here $\bar{\tilde{X}}$ is obtained from $\tilde{X}$ in the same way
that $\bar{X}$ was defined from $X$.

We remark that a similar version of Theorem~\ref{theoremsingular} has
appeared in \cite{TVcir} to deal with random matrices of i.i.d.
entries; see also \cite{BCC,NgV-comb} and the references therein.
However, our task here looks much harder as the entries in each row and
each column are not independent. We will now sketch the proof of
Theorem~\ref{theoremsingular}; more details will be presented in
Section~\ref{sectionstep1}.

Assume that $\sigma_{n-1}(\bar{X}+F)\le n^{-A}$. Then, by letting
$C=(c_{ij})_{2\le i,j\le n}$ be the cofactor matrix of $\bar{X}+F$,
there exist vectors $\Bx$ and $\By$ such that $\|\Bx\|=1$ and $\|\By
\|
\le n^{-A}$ and
\[
C\By=\det(\bar{X}+F)\Bx.
\]

So
\[
\|C\By\|=\bigl|\det(\bar{X}+F)\bigr|.
\]

Thus by the Cauchy--Schwarz inequality, with a loss of a factor of $n$
in probability and without loss of generality, we can assume that
%
%
\begin{equation}
\label{eqnsumsquare} \sum_{j=2}^n
|c_{2j}|^2 \ge n^{2A-1}\bigl|\det(\bar{X}+F)\bigr|^2.
\end{equation}

In what follows we fix the matrix $X_{(n-2)\times(n-1)}$ generated by
the last $(n-2)$ rows and the last $(n-1)$ columns of $X$
[equivalently, we fix the last $(n-2)$ rows of $\bar{X}$].

Let $s_2,\ldots,s_n$ be the column sums of $X_{(n-2)\times(n-1)}$. By
Theorem~\ref{theorembound}, the probability that all $x_{11},\ldots
,x_{1n},x_{21},\ldots,x_{2n}$ are greater than $n^{-2B-2}$ is bounded
from below by $1-O(n^{-B})$, in which case we have
%
%
\begin{eqnarray}
\label{eqnsi} s_i &\le&1-n^{-2B-2}\qquad \mbox{for all } i\ge2\quad
\mbox{and}
\nonumber
\\[-8pt]
\\[-8pt]
\nonumber
 0&\le& s_1:=(n-2)- (s_2+\cdots+s_n)
\le1-n^{-2B-2}.
\end{eqnarray}

Thus it is enough to justify Theorem~\ref{theoremsingular}
conditioning on this event.

Next, given a sequence $s_2,\ldots,s_n$ satisfying \eqref{eqnsi}, we
will choose $x_2:=x_{22},\ldots,x_n:=x_{2n}$ uniformly and,
respectively,
from the interval $[0,1-s_2],\ldots,\break   [0,1-s_n]$ such that
%
%
\begin{equation}
\label{eqnconstraint} s_1 \le x_2+\cdots+x_n
\le1.
\end{equation}

The upper bound guarantees that $x_1:=x_{21}=1-(x_2+\cdots+x_n) \ge0$,
while the lower bound ensures that $x_{11}=1-s_1-x_{21}=x_2+\cdots+x_n
-s_1 \ge0$.

We now express $\det(\bar{X}+F)$ as a linear form of its first row
$(x_2-x_1+f_{22},\ldots,x_n-x_1+f_{2n})$,
\[
\det(\bar{X}+F) = \sum_{2\le j \le n} c_{2j}(
\bar{X}+F) (x_j-x_1+f_{2j}).
\]

By using the fact that $x_1=1-\sum_{2\le j\le n} x_j$ we can rewrite
the above as
%
%
\begin{equation}
\label{eqndetshort} \det(\bar{X}+F)=\sum_{2\le j \le n}
\biggl(c_{2j}+\sum_{2\le i\le n} c_{2i}
\biggr)x_j+c,
\end{equation}
where $c$ is a constant depending on the $c_{2j}$'s and $f_{2j}$'s.

Observe that
\[
\sum_{2\le j\le n} \biggl|c_{2j}+\sum
_{2\le i\le n} c_{2i}\biggr|^2 = \sum
_{2\le j
\le n} |c_{2j}|^2 + (n+1)\biggl|\sum
_{2\le j \le n}c_{2j}\biggr|^2 \ge\sum
_{2\le
j\le n} |c_{2j}|^2.
\]

Thus, by increasing $A$ if needed, we obtain from \eqref{eqnsumsquare}
and \eqref{eqndetshort} the following:
\[
\biggl|\sum_{2\le j \le n}x_ja_j+c\biggr| \le
n^{-A},
\]
where
%
%
\begin{equation}
\label{eqna} a_j:=\frac{c_{2j}+\sum_{2\le i\le n} c_{2i}}{(\sum_{2\le j\le
n}|c_{2j}+\sum_{2\le i\le n} c_{2i}|^2)^{1/2}}.
\end{equation}

Roughly speaking, our approach to prove Theorem~\ref{theoremsingular}
consists of two main steps:

\begin{itemize}
\item\textit{Inverse step}. Given the matrix $X_{(n-2)\times(n-1)}$
for which all the column sums $s_i$ satisfy \eqref{eqnsi}, assume that
\[
\P_{x_{2},\ldots,x_{n}} \biggl(\biggl|\sum_{2\le j\le n}
a_j x_j+c\biggr|\le n^{-A} \biggr)\ge
n^{-B},
\]
where the probability is taken over all $x_i,2\le i$ which
satisfy \eqref{eqnconstraint}. Then there is a strong structure among
the cofactors $c_{2j}$ of $X_{(n-2)\times(n-1)}$.

\item\textit{Counting step}. With respect to $X_{(n-2)\times(n-1)}$,
the probability that there is a strong structure among the cofactors
$c_{2j}$ is negligible.
\end{itemize}

We pause to discuss the structure mentioned in the inverse step. A set
$Q\subset\C$ is a \emph{GAP of rank $r$} if it can be expressed as in
the form
\[
Q= \bigl\{g_0+ k_1g_1 +
\cdots+k_r g_r| k_i\in\Z, K_i \le
k_i \le K_i' \mbox{ for all } 1 \leq i
\leq r\bigr\}
\]
for some $(g_0,\ldots,g_r) \in\C^{r+1}$ and $(K_1,\ldots,K_r),
(K'_1,\ldots,K'_r) \in\Z^r$.

It is convenient to think of $Q$ as the image of an integer box $B:= \{
(k_1, \ldots, k_r) \in\Z^r| K_i \le k_i \le K_i' \} $ under the linear
map $\Phi\dvtx  (k_1,\ldots, k_r) \mapsto g_0+ k_1g_1 + \cdots+ k_r g_r$.

The numbers $g_i$ are the \emph{generators} of $Q$, the numbers $K_i'$
and $K_i$ are the \emph{dimensions} of $Q$, and $\Vol(Q):= |B|$ is the
\emph{size} of $B$. We say that $Q$ is \emph{proper} if this map is one
to one, or equivalently if $|Q| = \Vol(Q)$. For nonproper GAPs, we of
course have $|Q| < \Vol(Q)$. If $-K_i=K_i'$ for all $i\ge1$ and
$g_0=0$, we say that $Q$ is \textit{symmetric}.

We are now ready to state both of our steps in details.

\begin{theorem}[(Inverse step)]\label{theoremstep1}
Let $0<\ep<1$ and $B>0$ be given constants. Assume that
\[
\P_{x_{2},\ldots,x_{n}} \biggl(\biggl|\sum_{2\le j\le n}
a_j x_j+c\biggr|\le n^{-A} \biggr)\ge
n^{-B}
\]
for some sufficiently large integer $A$, where $a_j$ are defined in
\eqref{eqna}, and $x_j$ are chosen uniformly from the intervals
$[0,1-s_i]$ such that the constraint \eqref{eqnconstraint} holds. Then
there exists a vector $\Bu=(u_2,\ldots,u_{n})\in\C^{n-1}$ which
satisfies the following properties:

\begin{itemize}
\item$1/2 \le\|\Bu\| \le2$ and $|\langle\Bu, \row_i(\bar{X}+F)
\rangle| \le n^{-A+\gamma+2}$ for all but the first row $\row_1(\bar
{X}+F)$ of $\bar{X}+F$.

\item All but $n'$ components $u_i$ belong to a GAP $Q$ (not
necessarily symmetric) of rank $r=O_{B,\ep}(1)$, and of cardinality
$|Q|=n^{O_{B,\ep}(1)}$.

\item All the real and imaginary parts of $u_i$ and of the generators
of $Q$ are rational numbers of the form $p/q$, where $|p|,|q| \le n^{2A+3/2}$.
\end{itemize}
\end{theorem}

In the second step of the approach we show that the probability for
$X_{(n-2)\times(n-1)}$ having the above properties is negligible.

\begin{theorem}[(Counting step)]\label{theoremstep2}
With respect to $X_{(n-2)\times(n-1)}$, or equivalently, with respect
to the last $(n-2)$ rows of $\bar{X}$, the probability that there
exists a vector $\Bu$ as in Theorem~\ref{theoremstep1} is $\exp
(-\Theta(n))$.
\end{theorem}

\begin{pf}
First, we show that the
number of structural vectors $\Bu$ described in Theorem~\ref
{theoremstep1} is bounded by $n^{O_{B,\ep}(n)+O_A(n^\ep)}$. Indeed,
because each GAP is determined by its generators and its dimensions,
and because all the real and complex parts of the genrators are of the
form $p/q$ where $|p|,|q| \le n^{2A+3/2}$, there are $n^{O_{A,B,\ep
}(1)}$ GAPs which have rank $O_{B,\ep}(1)$ and size $n^{O_{B,\ep}(1)}$.
Next, for each determined GAP $Q$ of size $n^{O_{B,\ep}(1)}$, there are
$|Q|^n=n^{O_{B,\ep}(n)}$ ways to choose the $u_i$ as its elements. For
the remaining $O(n^{\ep})$ exceptional $u_i$ that may not belong to
$Q$, there are $n^{O_{A}(n^\ep)}$ ways to choose them as numbers of the
form $p/q$ where $|p|,|q| \le n^{2A+3/2}$. Putting these together we
obtain the bound as claimed.

Second, as for each fixed structural vector $\Bu$ from Theorem~\ref
{theoremstep1} we have $|\langle\Bu, \row_i(\bar{X}+F) \rangle|
=O(n^{-A+\gamma+2})$ for all $2\le i\le n-1$. So
%
%
\begin{eqnarray}
\label{eqnX}\qquad \biggl|\sum_{2\le j} u_j(x_{ij}-x_{i1}+f_{ij})
\biggr|&=& \biggl|\sum_{2\le j} x_{ij}\biggl(u_j+
\sum_{2\le k}u_k\biggr)-\sum
_{2\le j}u_j+\sum_{2\le j}
u_j f_{ij} \biggr|
\nonumber
\\[-8pt]
\\[-8pt]
\nonumber
&=&O\bigl(n^{-A+\gamma+2}\bigr).
\end{eqnarray}

We next view this inequality for the matrix model $Y$ and $\bar{Y}$,
where $Y$ was introduced in Section~\ref{sectionproperties}, and
$\bar
{Y}$ is obtained from $Y$ in the same way as how $\bar{X}$ was defined
from $X$,
%
%
\begin{equation}
\label{eqnY} \biggl|\sum_{2\le j} \frac{1}{n}y_{ij}
\biggl(u_j+\sum_{2\le k}u_k
\biggr)-\sum_{2\le
j}u_j+\sum
_{2\le j} u_j f_{ij}\biggr|=O
\bigl(n^{-A+\gamma+2}\bigr).
\end{equation}
Observe that
\[
\sum_{2\le j \le k}\biggl|u_j+\sum
_{2\le k\le n}u_k\biggr|^2 \ge\sum
_{2\le k\le n} u_k^2 \ge1/4.
\]

Thus there exits $j_0$ such that
\[
\biggl|u_{j_0}+\sum_{2\le k\le n}u_k\biggr|
\ge1/2\sqrt{n}.
\]

It then follows that for each $i$, with room to spare,
\begin{eqnarray*}
&&\hspace*{-4pt} \P \biggl( \biggl|\sum_{2\le j} \frac{1}{n}y_{ij}
\biggl(u_j+\sum_{2\le
k}u_k
\biggr)-\sum_{2\le j}u_j+\sum
_{2\le j} u_j f_{ij} \biggr|=O
\bigl(n^{-A+\gamma
+2}\bigr) \biggr)
\\
&&\hspace*{-4pt}\qquad=\P_{y_{ij},j\neq j_0} \P_{y_{ij_0}} \biggl( \biggl|\frac
{1}{n}y_{ij_0}
\biggl(u_{j_0}+\sum_{2\le k\le n}u_k
\biggr)\\
&&\hspace*{98pt}{}+ \sum_{j\neq j_0} \frac
{1}{n}y_{ij}
\biggl(u_j+\sum_{2\le k\le n}u_k
\biggr)- \cdots\biggr |\\
&&\hspace*{126pt}\qquad=O\bigl(n^{-A+\gamma
+2}\bigr)|y_{ij,j\neq j_0} \biggr)
\\
&&\hspace*{-4pt}\qquad=O\bigl(n^{-A+\gamma+10}\bigr),
\end{eqnarray*}
where in the last conditional probability estimate we used the fact
that $y_{ij}$ are i.i.d. exponentials of mean one.

Hence, for each fixed structural vector $\Bu$, the probability $\P
_{\Bu
}$ that \eqref{eqnY} holds for all rows $\row_i(\bar{Y}+F), 2\le
i\le
n-1$, is bounded by
\[
\P_{\Bu}\le n^{(-A+\gamma+10)(n-2)}.
\]

Summing over structural vectors $\Bu$, we thus obtain the following
upper bound for the probability that there exists a structural vector
$\Bu$ for which \eqref{eqnY} holds for all rows $\row_i(\bar{Y}+F),
2\le i\le n-1$
\[
\sum_{\Bu}\P_{\Bu}\le
n^{O_{B,\ep}(n)+O_A(n^\ep)}n^{(-A+\gamma
+10)(n-2)}=O\bigl(n^{-An/2}\bigr),
\]
provided that $A$ is large enough.

To conclude the proof of Theorem~\ref{theoremstep2}, we use Lemma~\ref
{lemmarelation1} to pass from $Y$ and $\bar{Y}$ back to $X$ and
$\bar
{X}$. The probability that there exists a structural vector $\Bu$ for
which \eqref{eqnX} holds for all rows $\row_i(\bar{X}+F),2\le i\le
n-1$, is then bounded by $O(n^{-An/2 +4n})=O (\exp(-\Theta
(n))
)$, provided that $A$ is sufficiently large.
\end{pf}


\section{\texorpdfstring{Proof of Theorem \protect\ref{theoremstep1}}
{Proof of Theorem 3.3}}\label{sectionstep1}

We recall from the assumptions of Theorem~\ref{theoremstep1} that
%
%
\begin{equation}
\label{eqnconcentrationx} \P_{x_2,\ldots,x_{n}} \biggl(\biggl|\sum
_{j\ge2} a_j x_j+c\biggr|\le
n^{-A} \biggr)\ge n^{-B},
\end{equation}
where $x_2,\ldots,x_n$ are uniformly sampled from the interval
$[0,1-s_2],\ldots,[0,1-s_n]$, respectively, so that \eqref
{eqnconstraint} holds.

This is a large concentration inequality for linear forms of mildly
dependent random variables. Our first goal is to relax these dependencies.

\subsection{A simple reduction step} Let $E_n$ be the set of all
$(x_2,\ldots,x_n)$ uniformly sampled from $[0,1-s_2] \times\cdots
\times
[0,1-s_n]$ so that \eqref{eqnconstraint} holds. We recall from \eqref
{eqnsi} that $s_i\le1-n^{-2B-2}$.

Consider the event $s_1\le x'_2+\cdots+x'_{n} \le1$, where $x'_i$ are
independently and uniformly sampled from the interval $[0,1-s_i]$, respectively.

Note that $\E(x'_{2}+\cdots+ x'_{n}) = \sum_{2\le i\le
n}(1-s_i)/2=(1-s_1)/2$. Since the random variables $x'_i- (1-s_i)/2$
are symmetric and uniform, the density function $f(x)$ of $x'_2+\cdots
+x'_{n}$ is maximized at $(1-s_1)/2$ and decreases as $|x-(1-s_1)/2)|$
increases. Thus we have
\begin{eqnarray*}
\P \bigl(\bigl(x'_2,\ldots,x'_n
\bigr)\in E_n \bigr) &=& \P\bigl(s_1\le
x'_2+\cdots +x'_{n} \le 1
\bigr)
\\
&= &\int_{s_1}^1 f(x)\,dx = \frac{\int_{s_1}^1 f(x)\,dx}{\int_{0}^{(1-s_2)+\cdots+(1-s_n)} f(x)\,dx}
\\
&\ge&\frac{1-s_1}{(1-s_2)+\cdots+(1-s_n)} = \frac{1-s_1}{1+s_1}
\\
&=&\Omega\bigl(n^{-2B-2}\bigr),
\end{eqnarray*}
where we noted from \eqref{eqnsi} that $s_1\le1-n^{-2B-2}$.

Observe that if we condition on $s_n\le x'_2+\cdots+x'_{n}\le1$, then
the distribution of $(x'_2,\ldots,x'_n)$ is uniform over the set $E_n$.
It thus follows from \eqref{eqnconcentrationx} that
%
%
\begin{equation}
\label{eqnconcentration1} \P_{x'_2,\ldots,x'_{n}} \biggl(\biggl|\sum
_{j\ge2} a_j x'_j+c\biggr|\le
n^{-A} \biggr)\ge n^{-3B-2}.
\end{equation}

In the next step of the reduction, we divide the intervals $[0,1-s_i]$
into disjoint intervals $I_{i1},\ldots,I_{ik_i}$ of length $n^{-3B-2}$,
where $k_i=(1-s_i)/n^{-3B-2}$ (without loss of generality, we assume
that $k_i$ are integers). Next, to sample $x'_i$ uniformly from the
interval $[0,1-s_i]$ we first choose at random an interval from $\{
I_{i1},\ldots,I_{ik_i}\}$ and then sample $x'_i$ from it. In this way,
\eqref{eqnconcentration1} implies that there exist intervals
$I_{ij_i}, 2\le i\le n$, such that if $x'_i$ are chosen uniformly from
$I_{ij_i}$ then
%
%
\begin{equation}
\label{eqnconcentration2} \P_{x'_2,\ldots,x'_{n}} \biggl(\biggl|\sum
_{j\ge2} a_j x'_j+c\biggr|\le
n^{-A} \biggr)\ge n^{-3B-2}.
\end{equation}

Observe furthermore that, by shifting $c$ if needed, we can assume that
$I_{ij_i}=[0,n^{-3B-2}]$ for all $i$. Finally, by passing to
$x_i'':=n^{3B+2}x'_i$ and by decreasing $A$ to $A-(3B+2)$, we can
assume that all $x_i'$ are uniformly sampled from the interval $[0,1]$.

\subsection{High concentration of linear form} A classical result of
Erd\H{o}s \cite{E} and Littlewood--Offord \cite{LO} asserts that if
$b_i$ are real numbers of magnitude $|b_i|\ge1$, then the probability
that the random sum $\sum_{i=1}^n b_ix_i$ concentrates on an interval
of length one is of order $O(n^{-1/2})$, where $x_i$ are i.i.d. copies
of a Bernoulli random variable. This remarkable inequality has
generated an impressive amount of research, particularly from the early
1960s to the late 1980s. We refer the reader to \cite{H,Kle} and the
references therein for these developments.

Motivated by inverse theorems from additive combinatorics, Tao and Vu
studied the underlying reason as to why the concentration probability
of $\sum_{i=1}^n b_ix_i$ on a short interval is large. A closer look at
the definition of GAPs defined in the previous section reveals that if
$b_i$ are very \textit{close} to the elements of a $\mathit{GAP}$ of rank $O(1)$
and size $n^{O(1)}$, then the concentration probability of $\sum_{i=1}^n b_ix_i$ on a short interval is of order $n^{-O(1)}$, where
$x_i$ are i.i.d. copies of a Bernoulli random variable.

It has been shown by Tao and Vu \cite{TVbull,TVcir,TVcomp} in an
implicit way, and by the current author and Vu \cite{NgV} in a more
explicit way that these are essentially the only examples that have
high concentration probability.

We say that a complex number $a$ is $\delta$-close to a set $Q\subset
\C
$ if there exists $q\in Q$ such that $|a-q|\le\delta$.

\begin{theorem}[(Inverse Littlewood--Offord theorem for linear forms \cite{NgV}, Corollary~2.10)]
\label{theoremILOlinear} Let $0 <\ep< 1$
and $C>0$. Let
$ \beta>0$ be an arbitrary real number that may depend on $n$. Suppose
that $b_i=(b_{i1},b_{i2})$ are complex numbers such that $\sum_{i=1}^n
\|b_i\|^2=1$, and
\[
\sup_a\P_{\Bx} \Biggl(\biggl|\sum
_{i=1}^n b_ix_i-a\biggr| \le
\beta\Biggr)=\rho\ge n^{-C},
\]
where $\Bx=(x_1,\ldots,x_n)$, and $x_i$ are i.i.d. copies of random
variable $\xi$ satisfying $\P(c_1\le\xi-\xi' \le c_2) \ge c_3$ for
some positive constants $c_1,c_2$ and $c_3$. Then, for any number $n'$
between $n^\ep$ and $n$, there exists a proper symmetric GAP $Q=\{\sum_{i=1}^r k_ig_i \dvtx k_i\in\Z, |k_i|\le L_i \}$ such that:

\begin{itemize}
\item at least $n-n'$ numbers $b_i$ are $\beta$-close to $Q$;

\item$Q$ has small rank, $r=O_{C,\ep}(1)$, and small cardinality
\[
|Q| \le\max \biggl(O_{C,\ep}\biggl(\frac{\rho^{-1}}{\sqrt{n'}}\biggr),1 \biggr);
\]

\item there exists a nonzero integer $p=O_{C,\ep}(\sqrt{n'})$ such
that all generators $g_i=(g_{i1},g_{i2})$ of $Q$ have the form
$g_{ij}=\beta\frac{p_{ij}} {p} $, with $p_{ij} \in\Z$ and
$|p_{ij}|= O_{C,\ep}(\beta^{-1} \sqrt{n'})$.
\end{itemize}
\end{theorem}

Theorem~\ref{theoremILOlinear} was proved in \cite{NgV} with
$c_1=1,c_2=2$ and $c_3=1/2$, but the proof there automatically extends
to any constants $0<c_1<c_2$ and $0<c_3$.

The interested reader is invited to also read \cite{Ng,RV,V} for other variants and further developments of such inverse results.

We now prove Theorem~\ref{theoremstep1}. Theorem~\ref
{theoremILOlinear} applied to \eqref{eqnconcentration2}, with
$n'=n^\ep,C=3B+2$ and $x_i$ being independently and uniformly
distributed over the interval $[0,1]$, implies that there exists a
vector $\Bv=(v_2,\ldots,v_n)$ such that:

\begin{itemize}
\item$|a_i-v_i|\le n^{-A}$ for all indices $i$ from $\{2,\ldots,n\}$;

\item all but $n'$ numbers $v_i$ belong to a GAP $Q$ of small rank,
$r=O_{B,\ep}(1)$, and of small cardinality $|Q|=O(n^{O_{B,\ep}(1)})$;

\item all the real and imaginary parts of $v_i$ and of the generators
of $Q$ are rational numbers of the form $p/q$, with $p,q \in\Z$ and
$|p|,|q|=O_{B,\ep}(n^{A+1/2})$.
\end{itemize}

Recall that
\[
a_j=\frac{c_{2j}+\sum_{2\le i\le n} c_{2i}}{(\sum_{2\le j\le
n}|c_{2j}+\sum_{2\le i\le n} c_{2i}|^2)^{1/2}}.
\]

We will translate the above useful information on the $a_j$'s to the
$c_j$'s. To do so we fist find a number of the form $p/n^A$, where
$p\in\Z$ and $-n^{A}\le p \le n^A$ such that
\[
\biggl|\frac{p}{n^{A}} -\frac{\sum_{2\le j\le n} c_{2j}}{(\sum_j|c_{2j}+\sum_{2\le i\le n} c_{2i}|^2)^{1/2}} \biggr|\le\frac{1}{n^{A}}.
\]

Thus, by shifting the GAP $Q$ by $p/n^A$, we obtain $|a_j'-v_j'|\le
2n^{-A}$, and so
\[
\bigl\|\Ba'-\Bv'\bigr\|=O\bigl(n^{-A+1/2}\bigr),
\]
where $\Ba'=(a_2',\ldots,a_n'), \Bv'=(v_2',\ldots,v_n')$ and
\[
a_j'=\frac{c_{2j}}{(\sum_j|c_{2j}+\sum_{2\le i\le n} c_{2i}|^2)^{1/2}}\quad \mbox{as well as}\quad
v_j'=v_j-\frac{p}{n^{A}}.
\]

By definition, $1/2n^2 \le\sum|a_j'|^2\le1$, so by the triangle inequality
\[
\bigl\|\Bv'\bigr\|\ge\bigl\|\Ba'\bigr\|-O\bigl(n^{-A+1/2}\bigr)
\ge1/\sqrt{2}n - O\bigl(n^{-A+1/2}\bigr)
\]
and
\[
\bigl\|\Bv'\bigr\|\le\bigl\|\Ba'\bigr\|+O\bigl(n^{-A+1/2}\bigr)
\le1+O\bigl(n^{-A+1/2}\bigr).
\]

More importantly, as $\Ba'$ is proportional to $(c_{22},\ldots,c_{2n})$
(which are the cofactors of $\bar{X}+F$), $\Ba'$ is orthogonal to all
but the first row of $\bar{X}+F$. In other words, $|\langle\Ba',\row
_i(\bar{X}+F) \rangle|=0$ for all $i\ge2$. It is thus implied that
\[
\bigl|\bigl\langle\Bv',\row_i(\bar{X}+F) \bigr\rangle\bigr|\le
n^{-A+\gamma+1}.
\]

In the last step of the proof, we find nonzero numbers $p',q'\in\Z,
|p'|,|q'|=O(n)$ so that $\|\Bv'\|/2 \le p'/q' \le2 \|\Bv'\|$.

Set
\[
\Bu:=\frac{q'}{p'} \Bv',
\]
and we then have:

\begin{itemize}
\item$1/2 \le\|\Bu\| \le2$ and $\langle\Bu, \row_i(\bar{X}+F)
\rangle\le n^{-A+\gamma+2}$ for all but the first rows of $\bar{X}+F$;

\item all but $n'$ components $u_i$ belong to a GAP $Q'$ (not
necessarily symmetric) of small rank, $r=O_{B,\ep}(1)$, and of small
cardinality $|Q'|=O(n^{O_{B,\ep}(1)})$;

\item all the real and imaginary parts of $u_i$ and of the generators
of $Q'$ are rational numbers of the form $p/q$, with $p,q\in\Z$ and
$|p|,|q|=O_{B,\ep}(n^{2A+3/2})$.
\end{itemize}


\section{Spectral concentration of i.i.d. random covariance
matrices}\label{sectionconcentration}

From now on we will mainly focus on the bounded model $\tilde{X}$
rather than on $X$. This is the model where we can relate to $\tilde
{Y}$, a matrix of bounded i.i.d. entries (defined in Section~\ref{sectionproperties})
for which concentration results may easily apply.
Furthermore, by Corollary~\ref{corbound}, there is not much difference
between the two models $X$ and $\tilde{X}$.

Having learned from Corollary~\ref{corsingulartilde{X}} that $|\det
(\sqrt{n}\bar{\tilde{X}}-z_0I_{n-1})|$ is bounded away from zero, we
will show that $\frac{1}{n}\log|\det(\sqrt{n}\bar{\tilde
{X}}-z_0I_{n-1})|$ is well concentrated around its mean. This result
will then immediately imply Theorem~\ref{theoremconvergence}.

In order to study the concentration of $\det(\sqrt{n}\bar{\tilde
{X}}-z_0I_{n-1})$, we might first relate it to the counterpart $\bar
{\tilde{Y}}$. However, the entries of the later model are not
independent, and so certain well-known concentration results for i.i.d.
matrices are not applicable. To avoid this technical issue, we will
modify $\sqrt{n}\bar{\tilde{X}}$ as follows. Observe that
%
%
\begin{equation}
\label{eqnequalitydet} \det(\sqrt{n}\bar{\tilde{X}}-z_0I_{n-1})=
\frac{1}{\sqrt{n}}\det (\sqrt {n}\tilde{X}_{(n-1) \times n}-F_{z_0}),
\end{equation}
where $F_{z_0}$ is the deterministic matrix obtained from $z_0I_{n-1}$
by attaching $(-\sqrt{n}, \ldots,-\sqrt{n})$ and $(-\sqrt{n}, 0,\ldots,
0)^T$ as its first row and first column, respectively, and $\tilde
{X}_{(n-1)\times n}$ is the matrix obtained from $\tilde{X}$ by
replacing its first row by a zero vector,
\[
\sqrt{n}\tilde{X}_{(n-1) \times n}-F_{z_0}:= \pmatrix{ \sqrt{n} &
\sqrt{n}& \cdots&\sqrt{n}\vspace*{2pt}
\cr
\sqrt{n} \tilde{x}_{21} &
\sqrt{n}\tilde{x}_{22}-z_0 & \cdots& \sqrt {n}
\tilde{x}_{2n} \vspace*{2pt}
\cr
\vdots& \vdots& \ddots& \vdots
\vspace*{2pt}
\cr
\sqrt{n} \tilde{x}_{n1} &\sqrt{n}\tilde{x}_{n2}
& \cdots& \sqrt {n}\tilde{x}_{nn}-z_0 }.
\]

As it turns out, it is more pleasant to work with $\tilde
{X}_{(n-1)\times n}$ because the entries of its counterpart $\tilde
{Y}_{(n-1)\times n}$ are now independent. To relate the singularity of
$\sqrt{n}\bar{\tilde{X}}-z_0I_{n-1}$ to that of $ \sqrt{n}\tilde
{X}_{(n-1) \times n}-F_{z_0}$, we have a crucial observation below.

\begin{claim}\label{claimspectralrelation} Suppose that $A$ is a
sufficiently large constant. We have
\[
\sigma_n( \sqrt{n}\tilde{X}_{(n-1) \times n}-F_{z_0}) \ge
\frac
{1}{n}\min \biggl(\frac{1}{\sqrt{2n}}\sigma_{n-1}(\sqrt{n}\bar
{\tilde {X}}-z_0I_{n-1}) -O\bigl(n^{-A}
\bigr),n^{-A} \biggr).
\]
\end{claim}

To prove this claim, let $\col_1,\ldots,\col_n$ be the columns of $
\sqrt {n}\tilde{X}_{(n-1) \times n}-F_{z_0}$. Let $\Bv=(v_1,\ldots,v_n)$ be
any unit vector. If $|v_1+\cdots+v_n|\ge n^{-A-1/2}$, then it is clear
that $\| (\sqrt{n}\tilde{X}_{(n-1) \times n}-F_{z_0})\Bv\|\ge|\sqrt {n}(v_1+\cdots+v_n)|\ge n^{-A}$. Otherwise, as $|v_1|^2+\cdots
+|v_n|^2=1$, we can easily deduce that $|v_2|^2+\cdots+|v_n|^2 \ge
1/2n$. Next, by the triangle inequality,
\begin{eqnarray*}
&&\bigl\| (\sqrt{n}\tilde{X}_{(n-1) \times n}-F_{z_0})v\bigr\|\\
&&\qquad=\biggl\|\sum
_{2\le
i\le
n} v_i \col_i\biggr\|= \biggl\|\sum
_{2\le i \le n}v_i(\col_i-
\col_1)+(v_1+\cdots +v_n)\col_1\biggr\|
\\
&&\qquad\ge\biggl\|\sum_{2\le i\le n}v_i
\col_i\biggr\| -n^{-A-1/2}\|\col_1\|
\\
&&\qquad\ge\bigl(|v_2|^2+\cdots+|v_n|^2
\bigr)^{1/2}\sigma_{n-1}(\sqrt{n}\bar{\tilde
{X}}-z_0I_{n-1})-\sqrt{2}n^{-A}
\\
&&\qquad\ge\frac{1}{\sqrt{2n}}\sigma_{n-1}(\sqrt{n}\bar{\tilde
{X}}-z_0I_{n-1}) - O\bigl(n^{-A}\bigr).
\end{eqnarray*}

Claim~\ref{claimspectralrelation} guarantees that the polynomial
probability bound for\break  $\sigma_{n-1}(\sqrt{n}\bar{\tilde
{X}}-z_0I_{n-1})$ from Corollary~\ref{corsingulartilde{X}} continues
to hold for\break  $\sigma_n( \sqrt{n}\tilde{X}_{(n-1) \times n}-F_{z_0})$
(with probably a worse value of $A$).

\begin{theorem}\label{theoremsingularX} There exists a positive
constant $A$ such that
\[
\P\bigl(\sigma_n(\sqrt{n}\tilde{X}_{(n-1) \times n}-F_{z_0})
\le n^{-A}\bigr)=O\bigl(n^{-3}\bigr).
\]
\end{theorem}

Our goal is then to establish a large concentration of\break  $\frac
{1}{n}\log
|\det(\sqrt{n}\tilde{X}_{(n-1)\times n} -F_{z_0})|$ around its mean. We
now consider $\tilde{Y}$.

\subsection{\texorpdfstring{Large concentration for $\tilde{Y}$}{Large concentration for Y}} Consider the i.i.d.
matrices $\tilde{Y}$ defined from Section~\ref{sectionproperties}, and
let $\tilde{Y}_{(n-1)\times n}$ be the matrix obtained from $\tilde{Y}$
by replacing its first row by the zero vector.

We first observe from Claim~\ref{claimspectralrelation} that
\begin{eqnarray*}
&&\sigma_n\biggl( \frac{1}{\sqrt{n}}\tilde{Y}_{(n-1) \times n}-F_{z_0}
\biggr) \\
&&\qquad\ge \frac{1}{n}\min \biggl(\frac{1}{\sqrt{2n}}\sigma_{n-1}
\biggl(\frac
{1}{\sqrt {n}}\bar{\tilde{Y}}-z_0I_{n-1}\biggr)
-O\bigl(n^{-A}\bigr),n^{-A} \biggr),
\end{eqnarray*}
where
\[
\frac{1}{\sqrt{n}}\tilde{Y}_{(n-1) \times n}-F_{z_0}= \pmatrix{
\sqrt{n} &\sqrt{n}& \cdots&\sqrt{n}\vspace*{2pt}
\cr
\frac{1}{\sqrt{n}}
\tilde{y}_{21} &\frac{1}{\sqrt{n}}\tilde {y}_{22}-z_0
& \cdots& \frac{1}{\sqrt{n}}\tilde{y}_{2n} \vspace*{2pt}
\cr
\vdots&
\vdots& \ddots& \vdots\vspace*{2pt}
\cr
\frac{1}{\sqrt{n}} \tilde{y}_{n1}
&\frac{1}{\sqrt{n}}\tilde {y}_{n2} & \cdots& \frac{1}{\sqrt{n}}
\tilde{y}_{nn}-z_0 }.
\]

On the other hand, conditioning on $\tilde{y}_{21},\ldots,\tilde
{y}_{n1}$, the entries $\tilde{y}_{ij}-\tilde{y}_{i1}$ of the matrix
$\bar{\tilde{Y}}$ are independent, and so we can apply known
singularity bounds, for instance \cite{TVcir'}, Theorem~2.1, for i.i.d.
matrices to conclude that for any positive constant $B$, there exists a
positive constant $A$ such that $\P(\sigma_{n-1}(\frac{1}{\sqrt {n}}\bar
{\tilde{Y}}-z_0I_{n-1})\le n^{-A})=O(n^{-B})$. Returning to $\tilde
{Y}_{(n-1) \times n}$, we hence obtain the following.

\begin{theorem}\label{theoremsingularY} For any positive constant
$B$, there exists a positive constant~$A$ such that
\[
\P\biggl(\sigma_n\biggl( \frac{1}{\sqrt{n}}\tilde{Y}_{(n-1) \times
n}-F_{z_0}
\biggr)\le n^{-A}\biggr)=O\bigl(n^{-B}\bigr).
\]
\end{theorem}

This bound will be exploited later on.

Next, let $H$ denote the following Hermitian matrix:
\[
H:=\biggl(\frac{1}{\sqrt{n}}\tilde{Y}_{(n-1)\times n}-F_{z_0}
\biggr)^\ast\biggl(\frac
{1}{\sqrt{n}}\tilde{Y}_{(n-1)\times n}-F_{z_0}
\biggr).
\]

It is clear that the eigenvalues $\lambda_1(H),\ldots, \lambda_n(H)$ of
$H$ can be written as
\[
\lambda_1(H)=\sigma_1^2\biggl(
\frac{1}{\sqrt{n}}\tilde{Y}_{(n-1)\times
n}-F_{z_0}\biggr),\ldots,
\lambda_n(H)=\sigma_n^2\biggl(
\frac{1}{\sqrt{n}}\tilde {Y}_{(n-1)\times n}-F_{z_0}\biggr),
\]
where $\sigma_i(\frac{1}{\sqrt{n}}\tilde{Y}_{(n-1)\times n}-F_{z_0})$
are the singular values of $\frac{1}{\sqrt{n}}\tilde{Y}_{(n-1)\times n}
- F_{z_0}$.

The following concentration result will serve as our main lemma.

\begin{lemma}\label{lemmaconcentrationgeneral} Assume that $f$ is a
function so that $g(x):=f(x^2)$ is convex and has finite Lipshitz norm
$\|g\|_L$. Then for any $\delta\ge C K\|g\|_L/n $, where $K=10\log n$
is the upper bound for the entries of $\tilde{Y}_{(n-1)\times n}$ and
$C$ is a sufficiently large absolute constant, we have
\[
\P \Biggl( \Biggl|\sum_{i=1}^nf\bigl(
\lambda_i(H)\bigr)-\E\Biggl(\sum_{i=1}^nf
\bigl(\lambda _i(H)\bigr)\Biggr) \Biggr|\ge\delta n \Biggr) =O \biggl(\exp
\biggl(-C'\frac{n^2\delta^2}{K^2
\|g\|_L^2}\biggr) \biggr);
\]
here $C'$ and the implied constant depend on $C$.
\end{lemma}

We remark that when $F_{z_0}$ vanishes, Lemma~\ref
{lemmaconcentrationgeneral} is essentially \cite{GZ}, Corollary~1.8,
of Guionnet and Zeitouni. We will show that the method there can be
easily extended for any deterministic matrix $F_{z_0}$.

\begin{pf*}{Proof of Lemma~\ref{lemmaconcentrationgeneral}} Consider the
following Hermitan matrix $K_{2n}$ of size $2n\times2n$
\[
K_{2n} = \pmatrix{ 0 & \displaystyle\biggl(
\frac{1}{\sqrt{n}}\tilde{Y}_{(n-1)\times n}-F_{z_0}\biggr)^\ast
\vspace*{2pt}\cr
\displaystyle\frac{1}{\sqrt{n}}\tilde{Y}_{(n-1)\times n}-F_{z_0} & 0 }.
\]

Apparently,
\begin{eqnarray*}
K_{2n}^2 &=&\left(\matrix{
\displaystyle\biggl(\frac{1}{\sqrt{n}}\tilde{Y}_{(n-1)\times n}-F_{z_0}
\biggr)^\ast\biggl(\frac
{1}{\sqrt{n}}\tilde{Y}_{(n-1)\times n}-F_{z_0}
\biggr)\vspace*{2pt}\cr
 0}\right.\\
&&\hspace*{7pt} \left.\matrix{
0 \vspace*{2pt}\cr \displaystyle\biggl(\frac{1}{\sqrt{n}}\tilde{Y}_{(n-1)\times n}-F_{z_0}\biggr)
\biggl(\frac
{1}{\sqrt {n}}\tilde{Y}_{(n-1)\times n}-F_{z_0}
\biggr)^\ast }\right).
\end{eqnarray*}

So to prove Lemma~\ref{lemmaconcentrationgeneral}, it is enough to
show that
%
%
\begin{eqnarray}
\label{eqnT} &&\P \Biggl( \Biggr|\sum_{i=1}^{2n}g
\bigl(\lambda_i(K_{2n})\bigr)-\E\Biggl(\sum
_{i=1}^{2n}g\bigl(\lambda_i(K_{2n})
\bigr)\Biggr) \Biggr|\ge2\delta n \Biggr)
\nonumber
\\[-8pt]
\\[-8pt]
\nonumber
&&\qquad=O \biggl(\exp \biggl(-C'
\frac{n^2\delta^2}{K^2 \|g\|_L^2}\biggr) \biggr),
\end{eqnarray}
where $\lambda_i(K_{2n})$ are the eigenvalues of $K_{2n}$.

Next, by following \cite{GZ}, Lemma~1.2, we obtain the following.

\begin{lemma}\label{lemmacv_lip}
The function $M \mapsto\tr(g(\frac{1}{\sqrt{n}}M+F))$ of Hermitian
matrices $M=(m_{ij})_{1\le i,j \le n}$, where $F$ is a deterministic
Hermitian matrix whose entries may depend on $n$, is a:

\begin{itemize}
\item convex function;
\item Lipschitz function of constant bounded by $2\|g\|_L$.
\end{itemize}
\end{lemma}

We refer the reader to Appendix \ref{sectioncvLip} for a proof of
Lemma~\ref{lemmacv_lip}. To deduce \eqref{eqnT} from Lemma~\ref
{lemmacv_lip}, we apply the following well-known Talagrand
concentration inequality \cite{T}.

\begin{lemma}\label{lemmaTalagrand}
Let $\D$ be the disk $\{z\in\C, |z|\le K\}$. For every product
probability $\mu$ in $\D^N$, every convex function $F\dvtx \C^N \mapsto
\R$
of Lipschitz norm $\|F\|_L$, and every $r\ge0$,
\[
\P\bigl(\bigl|F-M(F)\bigr|\ge r\bigr)\le4 \exp\bigl(-r^2/16K^2\|F
\|_L^2\bigr),
\]
where $M(F)$ denotes the median of $F$.
\end{lemma}

Indeed, let $F$ be the function $\dvtx \tilde{Y}' \mapsto\tr
(g(K_{2n}))=\tr
(g(\frac{1}{\sqrt{n}}\tilde{Y}'+F'))$, where
\[
\tilde{Y}'= \pmatrix{ 0 &
\tilde{Y}_{(n-1)\times n}^\ast
\vspace*{2pt}\cr
\tilde{Y}_{(n-1)\times n} & 0 }
\]
and
\[
F'= \pmatrix{ 0 &
-F_{z_0}^\ast
\vspace*{2pt}\cr
-F_{z_0} & 0 }.
\]

Observe that the entries of $\tilde{Y}'$ are supported on $|x|\le
K=10\log n$. By Lemma~\ref{lemmacv_lip}, $F$ is a convex function with
Lipschitz constant bounded by $2\|g\|_L$. The conclusion \eqref{eqnT}
of Lemma~\ref{lemmaconcentrationgeneral} then follows by applying
Lem\-ma~\ref{lemmaTalagrand}.
\end{pf*}

In what follows we will apply Lemma~\ref{lemmaconcentrationgeneral}
for two functions: one gives an almost complete control on the large
spectra of $H$, and the other yields a good bound on the number of
small spectra of $H$. We will choose $c$ to be a sufficiently small
constant, and with room to spare we set
\[
\ep=\delta=\Theta\bigl(n^{-c}\bigr).
\]

\subsection{Concentration of large spectra for i.i.d. matrices}
Following \cite{CV}\break  and~\cite{FRZ}, we first apply Lemma~\ref
{lemmaconcentrationgeneral} to the cut-off function $f_\ep(x):=\break \log
(\max(\ep,x))$. Note that $f_\ep(x^2)$ has Lipschitz constant $2\ep
^{-1/2}$. Although the function is not convex, it is easy to write it
as a difference of two convex functions of Lipschitz constant $O(\ep
^{-1/2})$, and so Lemma~\ref{lemmaconcentrationgeneral} applies
because $\delta=\Theta(n^{-c}) \ge C\ep^{1/2}K /n$.

\begin{theorem}\label{theoremconcentrationlargespectralY} We have
\begin{eqnarray*}
&&\P \biggl( \biggl|\sum_{\sigma_i^2(({1}/{\sqrt{n}})\tilde
{Y}_{(n-1)\times
n}-F_{z_0})\in S_\ep}\log\sigma_i
\biggl(\frac{1}{\sqrt{n}}\tilde {Y}_{(n-1)\times n}-F_{z_0}\biggr)\\
&&\hspace*{129pt}{}-\E
\biggl(\sum_{\sigma_i^2(\cdots) \in
S_\ep
}\log\sigma_i(\cdots)
\biggr) \biggr|\ge\delta n \biggr)
\\
&&\qquad= O \bigl(\exp\bigl(-n^2\delta^2\ep/K^2
\bigr) \bigr) = O\bigl(\exp\bigl(-n\log^2n\bigr)\bigr),
\end{eqnarray*}
where $S_\ep:=\{x\in\R, x\ge\ep\}$.
\end{theorem}

For short, from now on we set
\[
h_{\ep,\tilde{Y}_{(n-1)\times n}}(z_0):=\frac{1}{n}\E \biggl(\sum
_{\sigma
_i^2(({1}/{\sqrt{n}})\tilde{Y}_{(n-1)\times n}-F_{z_0})\in S_\ep
}\log \sigma_i\biggl(\frac{1}{\sqrt{n}}
\tilde{Y}_{(n-1)\times n}-F_{z_0}\biggr) \biggr).
\]

Serving as the main term, $h_{\ep,\tilde{Y}_{(n-1)\times n}}(z_0)$ will
play a key role in our analysis. In the next subsection we apply Lemma~\ref{lemmaconcentrationgeneral} to another function $f$.

\subsection{Concentration of the number of small eigenvalues for i.i.d.
matrices} Let~$I$ be the interval $[0,\ep]$. We are going to show that
the number $N_I$ of the eigenvalues $\lambda_i(H)$ which belong to $I$
is small with very high probability.\vadjust{\goodbreak}

It is not hard to construct two functions $f_1,f_2$ such that
$(f_1-f_2)-\1_I$ is nonnegative and supported on an interval of length
$\ep/C$, and so that both of $g_1(x)=f_1(x^2)$ and $g_2(x)=f_2(x^2)$
are convex functions of Lipschitz constant $O(\ep^{-1/2})$. (E.g., one may construct $f_1(x), f_2(x)$ in such a way that the even
function $g_1(x)=f_1(x^2)$ is identical to 1 on the interval $[-\ep
^{1/2},\ep^{1/2}]$ and being straight concave down from both edges with
a slope of $O(\ep^{-1/2})$, while the graph of the function
$g_2(x)=f_2(x^2)$ is obtained from that of $g_1(x)$ by replacing its
positive part with zero).

Next, by Lemma~\ref{lemmaconcentrationgeneral} we have
\[
\P \biggl( \biggl|\sum_{\lambda_i(H)}f_1\bigl(
\lambda_i(H)\bigr)-\E\biggl( \sum_{\lambda
_i(H)}f_1
\bigl(\lambda_i(H)\bigr)\biggr) \biggr| \ge\delta n \biggr)= O \bigl(\exp
\bigl(-n\log ^2n\bigr) \bigr)
\]
and
\[
\P\biggl ( \biggl|\sum_{\lambda_i(H)}f_2\bigl(
\lambda_i(H)\bigr)-\E\biggl( \sum_{\lambda
_i(H)}f_2
\bigl(\lambda_i(H)\bigr)\biggl)\biggr | \ge\delta n \biggr)= O \bigl(\exp \bigl(-n
\log ^2n\bigr) \bigr).
\]

By the triangle inequality, we thus have
\begin{eqnarray*}
&&\P \biggl( \biggl|\sum_{\lambda_i(H)}(f_1-f_2)
\bigl(\lambda_i(H)\bigr)-\E\biggl( \sum
_{\lambda_i(H)}(f_1-f_2) \bigl(
\lambda_i(H)\bigr)\biggr) \biggr| \ge2\delta n \biggr)\\
&&\qquad= O \bigl(\exp
\bigl(-n\log^2n\bigr) \bigr).
\end{eqnarray*}
Because the error-function $f=(f_1-f_2)-\1_I$ is nonnegative, it
follows that with probability $1- O(\exp(-n\log^2n))$
\[
\sum_{\lambda_i(H)}\1_I\bigl(
\lambda_i(H)\bigr) + \sum_{\lambda
_i(H)}f\bigl(
\lambda _i(H)\bigr) \le\E \biggl( \sum_{\lambda_i(H)}(f_1-f_2)
\bigl(\lambda_i(H)\bigr) \biggr) + 2\delta n,
\]
and hence
\begin{eqnarray*}
N_I=\sum_{\lambda_i(H)}\1_I\bigl(
\lambda_i(H)\bigr) &\le&\E \biggl( \sum_{\lambda
_i(H)}(f_1-f_2)
\bigl(\lambda_i(H)\bigr) \biggr) + 2\delta n
\\
&\le&2\E \biggl(\sum_{\lambda_i(H)}\1_J\bigl(
\lambda_i(H)\bigr) \biggr) + 2\delta n
\\
&\le&2\E(N_J)+2\delta n,
\end{eqnarray*}
where $J$ is the interval $[0,\ep+\ep/C]$ and $N_J$ is the number of
eigenvalues of $H$ in $J$. (Strictly speaking, we have to set $J=[-\ep
/C,\ep+\ep/C]$. However, as $\lambda_i$ are nonnegative, we can omit
its negative interval.)



To exploit the above information furthermore, we apply a result saying
that $N_J$ has small expected value (see also \cite{TVcov}, Proposition~28 and the references therein).

\begin{lemma}\label{lemmaspectralmean} For all $J\subset\R$ with
$|J|\ge K^2\log^2n/n^{1/2}$, one has
\[
N_J\ll n|J|
\]
with probability $1-\exp(-\omega(\log n))$. In particular,
\[
\E(N_J) \le Cn|J|,
\]
where $C$ is a sufficiently large constant.
\end{lemma}

We remark that this result holds for any deterministic matrix $F_0$ in
the definition of $H$. We defer the proof of Lemma~\ref
{lemmaspectralmean} to Appendix \ref{sectionNJ}.

In summary, we have obtained the following result.

\begin{theorem}\label{theoremconcentrationsmallspectralY} With
probability $O(\exp(-n\log^2n))$, we have
\[
N_I \ge2C\ep n + 2\delta n,
\]
where $N_I$ is the number of $\sigma_i(\frac{1}{\sqrt{n}}\tilde
{Y}_{(n-1)\times n}-F_{z_0})$ such that $\sigma_i^2(\frac{1}{\sqrt {n}}\tilde{Y}_{(n-1)\times n}-F_{z_0}) \in[0,\ep]$.
\end{theorem}

Consequently, it follows from Theorems \ref{theoremsingularY} and
\ref
{theoremconcentrationsmallspectralY} that with probability
$1-O(n^{-B})$ the following holds:
\begin{eqnarray*}
\frac{1}{n}\sum_{\sigma_i^2(({1}/{\sqrt{n}})\tilde
{Y}_{(n-1)\times
n}-F_{z_0}) \in[0,\ep]}\log
\sigma_i\biggl(\frac{1}{\sqrt{n}}\tilde {Y}_{(n-1)\times n}-F_{z_0}
\biggr)&=&O\bigl((\ep+\delta)\log n\bigr)\\
&=&O\bigl(n^{-c}\log n\bigr).
\end{eqnarray*}

Thus, combining with Theorem~\ref
{theoremconcentrationlargespectralY}, we infer the following.

\begin{theorem}\label{theoremlogdetY}
Let $z_0$ be fixed, and let $B$ be a positive constant. Then the
following holds with probability $1-O(n^{-B})$:
\begin{eqnarray*}
\biggl\llvert \frac{1}{n}\log\biggl|\det\biggl(\frac{1}{\sqrt{n}}\tilde
{Y}_{(n-1)\times
n}-F_{z_0}\biggr)\biggr| - h_{\ep,\tilde{Y}_{(n-1)\times n}}(z_0)
\biggr\rrvert &\le& 2\delta+ O\bigl(n^{-c}\log n\bigr) \\
&=&O
\bigl(n^{-c}\log n\bigr),
\end{eqnarray*}
where the implied constants depend on $B$.
\end{theorem}

\subsection{\texorpdfstring{Asymptotic formula for $h_{\ep,\tilde{Y}_{(n-1)\times n}}(z_0)$}
{Asymptotic formula for h epsilon, Y(n-1) x n (z_0)}}

We next claim that\break  $\frac{1}{n}\log|\det(\frac{1}{\sqrt{n}}\tilde
{Y}_{(n-1)\times n}-F_{z_0})|$ also converges to the corresponding part
of the circular law, and so gives an asymptotic formula for $h_{\ep
,\tilde{Y}_{(n-1)\times n}}(z_0)$.

\begin{theorem}\label{theoremlogdetYcir}
For almost all $z_0$, the following holds with probability one:
%
%
\begin{equation}
\label{eqnapproxcir1} \frac{1}{n}\log\biggl|\det\biggl(\frac{1}{\sqrt{n}}
\tilde{Y}_{(n-1)\times
n}-F_{z_0}\biggr)\biggr|- \int_{\C}
\log|w-z_0|\,d\mu_{\cir}(w)=o(1).\vadjust{\goodbreak}
\end{equation}
\end{theorem}

Note that this result is more or less the circular law for random
matrices with i.i.d. entries. To prove it we simply rely on \cite{TVcir}.

\begin{pf*}{Proof of Theorem~\ref{theoremlogdetYcir}} We first pass to
$\bar
{\tilde{Y}}$
\[
\bar{\tilde{Y}}= \pmatrix{ \tilde{y}_{22}-\tilde{y}_{21} &
\cdots& \tilde{y}_{2n}-\tilde{y}_{21}\vspace*{2pt}
\cr
\tilde{y}_{32}-\tilde{y}_{31} & \cdots&
\tilde{y}_{3n}-\tilde{y}_{31}\vspace*{2pt}
\cr
\vdots&
\vdots& \vdots\vspace*{2pt}
\cr
\tilde{y}_{n2}-\tilde{y}_{n1}
&\cdots& \tilde{y}_{nn}-\tilde{y}_{n1} },
\]
where $\tilde{y}_{ij}$ are i.i.d. copies of $\tilde{y}$.

As
\[
\det\biggl(\frac{1}{\sqrt{n}}\tilde{Y}_{(n-1)\times n}-F_{z_0}\biggr)=
\sqrt {n}\det \biggl(\frac{1}{\sqrt{n}}\bar{\tilde{Y}}-z_0I_{n-1}
\biggr),
\]
it is enough to prove the claim for $\det(\frac{1}{\sqrt{n}}\bar
{\tilde
{Y}}-z_0I_{n-1})$.

View $\bar{\tilde{Y}}$ as a sum of the matrix $(\tilde
{y}_{ij})_{2\le
i,j\le n}$ and $R$, the $(n-1)\times(n-1)$ matrix formed by $(-\tilde
{y}_{i1},\ldots,-\tilde{y}_{i1})$ for $2\le i\le n$. Because $R$ has
rank one and the average square of its entries $\frac{1}{n-1}\sum_{i}\tilde{y}_{i1}^2$ is bounded almost surely (with respect to
$\tilde
{y}_{21},\ldots,\tilde{y}_{n1}$), \cite{TVcir}, Corollary~1.15, applied
to $\bar{\tilde{Y}}$ implies that the ESD of $\frac{1}{\sqrt{n}}
\bar
{\tilde{Y}}$ converges almost surely to the circular law.

Finally, thanks to \cite{TVcir}, Theorem~1.20, for almost all $z_0$ the
following holds with probability one:
\[
\frac{1}{n}\log \biggl|\det\biggl(\frac{1}{\sqrt{n}}\bar{\tilde
{Y}}-z_0I_{n-1}\biggr) \biggr|-\int_{\C}
\log|w-z_0|\,d\mu_{\cir}(w) =o(1).
\]
\upqed\end{pf*}

Theorems \ref{theoremlogdetY} and \ref{theoremlogdetYcir}
immediately imply that for almost all $z_0$,
%
%
\begin{equation}
\label{eqnapproxcir} h_{\ep,\tilde{Y}_{(n-1)\times n}}(z_0)- \int
_{\C} \log|w-z_0|\,d\mu_{\cir}(w)=o(1).
\end{equation}

By substituting \eqref{eqnapproxcir} back into Theorem~\ref
{theoremconcentrationlargespectralY}, we have
\begin{eqnarray}\label{eqnconcentrationlargespectralYcir}
&&\P \biggl(\biggl |\frac{1}{n} \sum_{\sigma_i^2(({1}/{\sqrt
{n}})\tilde
{Y}_{(n-1)\times n}-F_{z_0})\in S_\ep}\log
\sigma_i\biggl(\frac{1}{\sqrt {n}}\tilde{Y}_{(n-1)\times n}-F_{z_0}
\biggr)\nonumber\\
&&\hspace*{137pt}{}-\int_{\C} \log|w-z_0|\,d\mu
_{\cir
}(w) \biggr| \ge\delta+o(1) \biggr)
\\
&&\qquad =O\bigl(\exp\bigl(-n\log^2n
\bigr)\bigr).\nonumber
\end{eqnarray}

\section{\texorpdfstring{Large concentration for $\tilde{X}$, proof of Theorem \protect\ref{theoremconvergence}}
{Large concentration for X, proof of Theorem 1.10}}\label{sectionproof}

In this section we will
apply the transference principle of Lemma~\ref{lemmarelation2} to
pass the results of Section~\ref{sectionconcentration} back to
$\tilde
{X}$. Our treatment here is similar to \cite{CDS}, Section~4.\vadjust{\goodbreak}

By Lemma~\ref{lemmarelation2} and \eqref
{eqnconcentrationlargespectralYcir}, conditioning on $\tilde{Y}\in
\tilde{D}_n$ we have
\begin{eqnarray}\label{eqnconcentrationYcondition}
&&\P \biggl( \biggl|\frac{1}{n}\sum_{\sigma_i^2(({1}/{\sqrt
{n}})\tilde
{Y}_{(n-1)\times n}-F_{z_0})\in S_\ep} \log
\sigma_i\biggl(\frac{1}{\sqrt {n}}\tilde{Y}_{(n-1)\times n}-F_{z_0}
\biggr) \nonumber\\
&&\hspace*{133pt}{}-\int_{\C} \log|w-z_0|\,d\mu
_{\cir
}(w)\biggr | \ge\delta+ o(1)| \tilde{Y}\in\tilde{D}_n
\biggr)
\\
&&\qquad =O\bigl(n^{10n}\exp\bigl(-n
\log^2n\bigr)\bigr)=O\bigl(\exp\bigl(-n\log^2n/2\bigr)
\bigr).\nonumber
\end{eqnarray}

Next, for each $\tilde{Y} \in\tilde{D}_n$ we will compare the singular
values of $\frac{1}{\sqrt{n}}\tilde{Y}_{(n-1)\times n}-F_{z_0}$ with
those of $\sqrt{n}\tilde{X}_{(n-1)\times n}-F_{z_0}$, where $\tilde{X}$
is determined by $\Phi(\frac{1}{n}\tilde{Y})$, that is, $\tilde
{x}_{ij}=\frac{1}{n}\tilde{y}_{ij}$ for all $2\le i,j\le n$.

By definition, as $\tilde{Y} \in\tilde{D}_n$, we have $|\frac
{1}{n}\tilde{y}_{i1}-\tilde{x}_{i1}|\le n^{-4}$, and so the operator
norm of the difference matrix is bounded by
\[
\biggl\|\biggl(\frac{1}{\sqrt{n}}\tilde{Y}_{(n-1)\times n}-F_{z_0}\biggr)-(
\sqrt {n}\tilde{X}_{(n-1)\times n}-F_{z_0}) \biggr\|\le\frac{1}{n^2}.
\]

This leads to a similar bound for the singular values for every $i$
(see, e.g.,~\cite{HJ})
%
%
\begin{equation}
\label{eqnsingularvaluescomparison} \biggl|\sigma_i\biggl(
\frac{1}{\sqrt{n}}\tilde{Y}_{(n-1)\times
n}-F_{z_0}\biggr)-
\sigma_i(\sqrt{n}\tilde{X}_{(n-1)\times n}-F_{z_0})\biggr |\le
\frac{1}{n^2}.
\end{equation}

Notice furthermore that, conditioning on $\tilde{Y}\in\tilde{D}_n$,
$\Phi(\frac{1}{n}\tilde{Y})$ is uniformly distributed on the set
$\tilde
{S}_n$ of bounded doubly stochastic matrices $\tilde{X}$. Thus, with a
slight modification to $\ep$ by an amount of $n^{-2}$ [thus the order
of $\ep$ remains $\Theta(n^{-c})$], we obtain from \eqref
{eqnconcentrationYcondition} the following upper tail bound with
respect to $\tilde{X}$:
\begin{eqnarray*}
&&\P \biggl( \frac{1}{n}\sum_{\sigma_i^2(\sqrt{n}\tilde
{X}_{(n-1)\times
n}-F_{z_0})\in S_{\ep+n^{-2}}} \log
\sigma_i(\sqrt{n}\tilde {X}_{(n-1)\times n}-F_{z_0})\\
&&\hspace*{84pt}{}-\int
_{\C} \log|w-z_0|\,d\mu_{\cir
}(w)\ge
\delta+ o(1) \biggr)
\\
&&\qquad=O\bigl(\exp\bigl(-n\log^2n/2\bigr)\bigr).
\end{eqnarray*}

Also, we obtain a similar probability bound for the lower tail
\begin{eqnarray*}
&&\P \biggl( \frac{1}{n}\sum_{\sigma_i^2(\sqrt{n}\tilde
{X}_{(n-1)\times
n}-F_{z_0})\in S_{\ep-n^{-2}}} \log
\sigma_i(\sqrt{n}\tilde {X}_{(n-1)\times n}-F_{z_0})\\
&&\hspace*{68pt}{}-\int
_{\C} \log|w-z_0|\,d\mu_{\cir
}(w)\le -
\bigl(\delta+ o(1)\bigr) \biggr)
\\
&&\qquad=O\bigl(\exp\bigl(-n\log^2n/2\bigr)\bigr).
\end{eqnarray*}

Notice that these bounds hold for any $\ep=\Theta(n^{-c})$. By gluing
them together we infer the following variant of \eqref
{eqnconcentrationYcondition}.

\begin{theorem}\label{theoremconcentrationlargespectralX}
With respect to $\tilde{X}$ we have
\begin{eqnarray*}
&&\P \biggl( \biggl|\frac{1}{n}\sum_{\sigma_i^2(\sqrt{n}\tilde
{X}_{(n-1)\times
n}-F_{z_0})\in S_\ep} \log
\sigma_i(\sqrt{n}\tilde{X}_{(n-1)\times
n}-F_{z_0}) \\
&&\hspace*{116pt}{}-\int
_{\C} \log|w-z_0|\,d\mu_{\cir}(w) \biggr| \ge
\delta+ o(1) \biggr)
\\
&&\qquad=O\bigl(\exp\bigl(-n\log^2n/2\bigr)\bigr).
\end{eqnarray*}
\end{theorem}

Next, conditioning on $\tilde{Y}\in\tilde{D}_n$, by Theorem~\ref
{theoremconcentrationsmallspectralY} and Lemma~\ref
{lemmarelation2}, with probability $O(n^{10n}\exp(-n\log
^2n))=O(\exp
(-n\log^2n/2))$ we have
\[
N_I \ge2C\ep n + 2\delta n,
\]
where $N_I$ is the number of $\sigma_i(\frac{1}{\sqrt{n}}\tilde
{Y}_{(n-1)\times n}-F_{z_0})$ such that $\sigma_i^2(\frac{1}{\sqrt {n}}\tilde{Y}_{(n-1)\times n}-F_{z_0}) \in[0,\ep]$.

Because $\Phi(\frac{1}{n}\tilde{Y})$ is uniformly distributed on the
set $\tilde{S}_n$ conditioning on $\tilde{Y}\in\tilde{D}_n$, and also
because of \eqref{eqnsingularvaluescomparison}, we imply the following.

\begin{theorem}\label{theoremconcentrationsmallspectralX} With
probability $O(\exp(-n\log^2n))$ with respect to $\tilde{X}$, we have
\[
N_I \ge2C\biggl(\ep+\frac{1}{n^2}\biggr) n + 2\delta n,
\]
where $N_I$ is the number of $\sigma_i(\sqrt{n}\tilde
{X}_{(n-1)\times
n}-F_{z_0})$ such that $\sigma_i^2(\sqrt{n}\tilde{X}_{(n-1)\times
n}-F_{z_0}) \in[0,\ep]$.
\end{theorem}

We now gather the ingredients together to complete the proof of our
main result.

\begin{pf*}{Proof of Theorem~\ref{theoremconvergence} for $\tilde{X}$}
By Theorems \ref{theoremsingularX} and \ref
{theoremconcentrationsmallspectralX}, we have that
\begin{eqnarray*}
&&\P \biggl(\frac{1}{n}\sum_{\sigma_i^2(\sqrt{n}\tilde
{X}_{(n-1)\times
n}-F_{z_0}) \in[0,\ep]}\log
\sigma_i(\sqrt{n}\tilde {X}_{(n-1)\times
n}-F_{z_0})=O
\bigl((\ep+\delta)\log n\bigr) \biggr)\\
&&\qquad=1-O\bigl(n^{-3}\bigr).
\end{eqnarray*}

A combination of this fact with Theorem~\ref
{theoremconcentrationlargespectralX} implies that for almost all~$z_0$,
\begin{eqnarray*}
&&\P \biggl( \biggl|\frac{1}{n}\log|\det(\sqrt{n}\tilde {X}_{(n-1)\times
n}-F_{z_0})
- \int_{\C} \log|w-z_0|\,d\mu_{\cir}(w)
\biggr|=o(1) \biggr)\\
&&\qquad =1-O\bigl(n^{-3}\bigr).
\end{eqnarray*}

Hence, by \eqref{eqnequalitydet},
\[
\P \biggl( \biggl|\frac{1}{n}\log|\det(\sqrt{n}\bar{\tilde{X}}-z_0I_{n-1})
- \int_{\C} \log|w-z_0|\,d\mu_{\cir}(w)
\biggr|=o(1) \biggr) =1-O\bigl(n^{-3}\bigr),
\]
completing the proof.
\end{pf*}

\begin{appendix}\label{app}
\section{\texorpdfstring{Proof of Lemma~\lowercase{\protect\ref{lemmacv_lip}}}
{Proof of Lemma 5.6}}\label{sectioncvLip}
The main goal of this section is to justify Lemma~\ref{lemmacv_lip}.
Although our proof is identical to \cite{GZ}, Theorem~1.1 and \cite{GZ}, Corollary~1.8, let us present it here for the sake of completeness.

\subsection{Convexity} For simplicity, we first show that the function
$M\mapsto \tr(g(M+F))$ is convex. It then follows that the function $M
\mapsto\tr(g(\frac{1}{\sqrt{n}}M+F))$ is also convex.

For any Hermitian matrices $U$ and $V$
\[
g(V+F)-g(U+F)= \int_0^1 Dg \bigl(U+F+
\eta(V-U) \bigr)\sharp(V-U)\,d\eta
\]
where
\[
Dg(U+F)\sharp(V) = \lim_{\ep\rightarrow0} \ep^{-1} \bigl(g(U+F+
\ep V)-g(U+F) \bigr).
\]

For polynomial functions $g$, the noncommutative derivation $D$ can be
computed, and one finds in particular that for any $p\in\N$,
%
%
\begin{eqnarray}
\label{eqnp}&& (V+F)^p -(U+F)^p
\nonumber
\\
&&\qquad=\int_{0}^1
\Biggl( \sum_{k=0}^{p-1}\bigl(U+F+\eta (V-U)
\bigr)^k(V-U)\\
&&\hspace*{76pt}{}\times\bigl(U+F+\eta(V-U)\bigr)^{p-k-1} \Biggr) \,d
\eta.\nonumber
\end{eqnarray}

For such a polynomial function, by taking the trace and using $\tr
(AB)=\tr(BA)$, one deduces that
%
%
%
\begin{eqnarray}
\label{eqnp'}&& \tr \bigl((U+F)^p \bigr) - \tr \biggl(\biggl(
\frac{U+V}{2} +F\biggr)^p \biggr)
\nonumber
\\[-8pt]
\\[-8pt]
\nonumber
&&\qquad= p \int_{0}^1
\tr \biggl(\biggl(\frac{U+V}{2}+F+\eta\frac{U-V}{2}
\biggr)^{p-1}\frac{U-V}{2} \biggr) \,d\eta,
\\
\label{eqnp''} &&\tr \bigl((V+F)^p \bigr) - \tr \biggl(\biggl(
\frac{U+V}{2} +F\biggr)^p \biggr)
\nonumber
\\[-8pt]
\\[-8pt]
\nonumber
&&\qquad= p \int_{0}^1
\tr \biggl(\biggl(\frac{U+V}{2}+F-\eta\frac{U-V}{2}
\biggr)^{p-1}\frac{V-U}{2} \biggr) \,d\eta.
\end{eqnarray}

It follows from \eqref{eqnp}, \eqref{eqnp'} and \eqref{eqnp''} that
%
%
\begin{eqnarray}
\Delta&:=&\tr \bigl((U+F)^p \bigr)+ \tr \bigl((V+F)^p
\bigr) -2 \tr \biggl(\biggl(\frac
{U+V}{2}+F\biggr)^p \biggr)
\nonumber
\\[-8pt]
\\[-8pt]
\nonumber
&=&\frac{p}{2} \sum_{k=0}^{p-2} \int
_{0}^1 \int_0^1
\eta \,d \eta \,d\theta \tr \bigl((U-V)Z_{\eta,\theta}^k(U-V)Z_{\eta,\theta}^{p-2-k}
\bigr)
\end{eqnarray}
with
\[
Z_{\eta,\theta}:= \frac{U+V}{2}+F -\eta\frac{U-V}{2} + \eta
\theta(U-V).
\]

Next, for fixed $\eta,\theta\in[0,1]^2$, and fixed $U,V,F$ Hermitian
matrices, $Z_{\eta,\theta}$ is also Hermitian, and so we can find a
unitary matrix $U_{\eta,\theta}$ and a diagonal matrix $D_{\eta
,\theta
}$ with real diagonal entries $\lambda_{\eta,\theta}(1),\ldots,
\lambda
_{\eta,\theta}(n)$ so that
\[
Z_{\eta,\theta}=U_{\eta,\theta}D_{\eta,\theta}U^\ast_{\eta
,\theta}.
\]

Let $W_{\eta,\theta}= U_{\eta,\theta}= U^\ast_{\eta,\theta} (U-V)
U_{\eta,\theta}$. Then
%
%
\begin{eqnarray}
\label{eqnDelta} \Delta&=& \frac{p}{2} \sum_{k=0}^{p-2}
\int_{0}^1 \int_0^1
\eta \,d \eta \,d\theta\,\tr \bigl(W_{\eta,\theta} D_{\eta,\theta}^k
W_{\eta,\theta} D_{\eta,\theta}^{p-2-k} \bigr)
\nonumber
\\[-8pt]
\\[-8pt]
\nonumber
&=&\frac{p}{2} \sum_{k=0}^{p-2} \int
_{0}^1 \int_0^1
\eta \,d \eta \,d \theta \sum_{k=0}^{p-2} \sum
_{1\le i,j\le n} \lambda_{\eta,\theta
}^k(i)
\lambda _{\eta,\theta}^{p-2-k}(j) \bigl|W_{\eta,\theta}(ij)\bigr|^2.
\end{eqnarray}

But
\begin{eqnarray*}
\sum_{k=0}^{p-2}\lambda_{\eta,\theta}^k(i)
\lambda_{\eta,\theta
}^{p-2-k}(j) &=& \frac{\lambda_{\eta,\theta}^{p-1}(i) - \lambda
_{\eta
,\theta}^{p-1}(j)}{\lambda_{\eta,\theta}(i) -\lambda_{\eta,\theta}(j)} \\
&=& (p-1)\int
_{0}^1 \bigl(\alpha\lambda_{\eta,\theta}(j)+(1-
\alpha )\lambda _{\eta,\theta}(i) \bigr)^{p-2}\,d\alpha.
\end{eqnarray*}

Hence, substituting into \eqref{eqnDelta} gives
%
%
\begin{eqnarray}
\label{eqnx^p} \Delta&=& \frac{1}{2} \sum
_{1\le i,j\le n} \int_0^1 \int
_0^1 \int_0^1
\,d\alpha\eta \,d \eta \,d\theta\bigl|W_{\eta,\theta}(ij)\bigr|^2
\nonumber\\
&&\hspace*{89pt}{}\times g''\bigl(\alpha \lambda _{\eta,\theta}(j)+(1-
\alpha)\lambda_{\eta,\theta}(i)\bigr)\\
&\ge&0\nonumber
\end{eqnarray}
for the polynomial $g(x)=x^p$.\vadjust{\goodbreak}

Now, with $U,V,F$ being fixed, the eigenvalues $\lambda_{\eta,\theta
}(1),\ldots, \lambda_{\eta,\theta}(n)$ and the entries of $W_{\eta
,\theta
}$ are uniformly bounded. Hence, by Runge's theorem, we can deduce by
approximation that \eqref{eqnx^p} holds for any twice continuously
differentiable function $g$. As a consequence, for any such convex
function we have $g'' \ge0$ and
\[
\Delta= \tr \bigl(g(U+F) \bigr)+ \tr \bigl(g(V+F) \bigr) -2 \tr \biggl(g\biggl(
\frac
{U+V}{2}+F\biggr) \biggr)\ge0.
\]

\subsection{Boundedness} Now we show that the function $M \mapsto\tr
(g(\frac{1}{\sqrt{n}}M+F))$ has Lipschitz constant bounded by $2\|g\|_L$.

First, for any bounded continuously differentiable function $g$ we will
show that
\begin{eqnarray*}
&&\sum_{1\le i,j\le n}\!\biggl(\!d_{\Re(x_{ij})}\tr\biggl(g
\biggl(\frac{1}{\sqrt {n}}M+F\biggr)\!\biggr)\!\biggr)^2+ \sum
_{1\le i,j\le n}\!\biggl(\!d_{\Im(x_{ij})}\tr \biggl(g\biggl(
\frac
{1}{\sqrt{n}}M+F\biggr)\!\biggr)\!\biggr)^2 \\
&&\qquad\le4\|g
\|^2_L.
\end{eqnarray*}

We can verify that
%
%
\begin{equation}
\label{eqntraceformula} d_{\Re(x_{ij})}\tr \biggl(g\biggl(\frac{1}{\sqrt{n}}M+F
\biggr) \biggr) = \frac
{1}{\sqrt {n}}\tr \biggl(g'\biggl(
\frac{1}{\sqrt{n}}M+F\biggr)\Delta_{ij} \biggr),
\end{equation}
where $\Delta_{ij}(kl)=1$ if $kl=ij$ or $ji$ and zero otherwise.

Indeed, \eqref{eqntraceformula} is a consequence of \eqref{eqnp} for
polynomial functions, and it can be extended for bounded continuously
differentiable functions by approximations. In other words, we have
\begin{eqnarray*}
&& d_{\Re(x_{ij})}\tr \biggl(g\biggl(\frac{1}{\sqrt{n}}M+F\biggr) \biggr)\\
&&\qquad=
\cases{ \displaystyle\frac{1}{\sqrt{n}} \biggl(g'\biggl(\frac{1}{\sqrt{n}}M+F
\biggr) (ij)+ g'\biggl(\frac
{1}{\sqrt {n}}M+F\biggr) (ji) \biggr), &\quad
$i\neq j;$\vspace*{2pt}
\cr
\frac{1}{\sqrt{n}} g'\biggl(
\frac{1}{\sqrt{n}}M+F\biggr) (ii), & \quad $i=j.$}
\end{eqnarray*}

Hence,
\begin{eqnarray*}
\sum_{i,j} \biggl(d_{\Re(x_{ij})}\tr \biggl(g
\biggl(\frac{1}{\sqrt {n}}M+F\biggr) \biggr) \biggr)^2& \le&
\frac{2}{n} \sum_{i,j} \biggl|g'
\biggl(\frac{1}{\sqrt{n}}M+F\biggr) (ij)\biggr|^2 \\
&= &\frac
{2}{n} \tr
\biggl(g'\biggl(\frac{1}{\sqrt{n}}M+F\biggr) g'\biggl(
\frac{1}{\sqrt {n}}M+F\biggr)^\ast \biggr).
\end{eqnarray*}

But if $\lambda_1,\ldots,\lambda_n$ denote the eigenvalues of $\frac
{1}{\sqrt{n}}M+F$, then
\[
\tr \biggl(g'\biggl(\frac{1}{\sqrt{n}}M+F\biggr)g'
\biggl(\frac{1}{\sqrt{n}}M+F\biggr)^\ast \biggr) = \frac{1}{n}
\sum\bigl(g'(\lambda_i)\bigr)^2
\le\bigl\|g'\bigr\|_\infty^2.
\]
Thus we have
\[
\sum_{i,j} \biggl(d_{\Re(x_{ij})}\tr \biggl(g
\biggl(\frac{1}{\sqrt {n}}M+F\biggr) \biggr) \biggr)^2 \le2
\bigl\|g'\bigr\|_\infty^2.
\]

The same argument applies for derivatives with respect to $\Im
(x_{ij})$, and so by integration by parts and by the Cauchy--Schwarz inequality,
\[
\biggl|\tr \biggl(g\biggl(\frac{1}{\sqrt{n}}U+F\biggr) \biggr) -\tr \biggl(g\biggl(
\frac
{1}{\sqrt {n}}V+F\biggr) \biggr)\biggr | \le2\|g\|_L \|U-V\|
\]
for any $U$ and $V$.

Observe that the last result for bounded continuously differentiable
functions naturally extends to Lipschitz functions by approximation,
completing the proof.

\section{\texorpdfstring{Proof of Lemma \lowercase{\protect\ref{lemmaspectralmean}}}
{Proof of Lemma 5.11}}\label{sectionNJ}

Note that if $F_{z_0}$ vanishes, then Lemma~\ref{lemmaspectralmean}
is just \cite{TVcov}, Proposition~28; see also~\cite{BS}. We show that
the method there extends easily to any deterministic~$F_{z_0}$.

Assume for contradiction that
\[
|N_J| \ge C n |J|
\]
for some large constant $C$ to be chosen later. We will show that this
will lead to a contradiction with high probability.

We will control the eigenvalue counting function $N_J$ via the
Stieltjes transform
\[
s(z):=\frac{1}{n}\sum_{j=1}^n
\frac{1}{\lambda_j(H)-z}.
\]

Fix $J$ and let $x$ be the midpoint of $J$. Set $\eta:=|J|/2$ and
$z:=x+i\eta$, and we then have
\[
\Im\bigl(s(z)\bigr) \ge\frac{4}{5}\frac{N_J}{\eta n}.
\]

Hence,
%
%
\begin{equation}
\label{eqnim} \Im\bigl(s(z)\bigr) \gg C.
\end{equation}

Next, with $H':= (\frac{1}{\sqrt{n}}\Phi(\tilde{Y})-F_{z_0})(\frac
{1}{\sqrt{n}}\Phi(\tilde{Y})-F_{z_0})^\ast= \frac{1}{n}M M^\ast$ where
$M:=\Phi(\tilde{Y})-\sqrt{n}F_{z_0}$, we have (see also \cite{BS}, Chapter~11)
\[
s(z)= \frac{1}{n} \sum_{k\le n}
\frac{1}{h_{kk}' -z - \Ba_k^\ast
(H_k'-zI)^{-1}\Ba_k},
\]
where $h_{kk}'$ is the $kk$ entry of $H'$; $H_k'$ is the $n-1$ by $n-1$
matrix with the $k$th row and $k$th column of $H'$ removed; and $\Ba_k$
is the $k$th column of $H'$ with the $k$th entry removed.\vadjust{\goodbreak}

Note that $\Im(\frac{1}{z}) \le\frac{1}{\Im(z)}$, one concludes from
\eqref{eqnim} that
\[
\frac{1}{n} \sum_{k\le n} \frac{1}{ |\eta+\Im(\Ba_k^\ast
(H_k'-zI)^{-1} \Ba_k) |} \gg
C.
\]

By the pigeonhole principle, there exists $k$ such that
%
%
\begin{equation}
\label{eqnk} \frac{1}{ |\eta+ \Im(\Ba_k^\ast(H_k'-zI)^{-1}\Ba_k) |} \gg C.
\end{equation}

Fix such $k$, note that
\[
\Ba_k = \frac{1}{n}M_k\row_k^\ast\quad
\mbox{and} \quad H_k'=\frac
{1}{n}M_kM_k^\ast,
\]
where $\row_k=\row_k(M)$ and $M_k$ is the $(n-1)\times n$ matrix formed
by removing $\row_k(M)$ from $M$. Thus if we let $\Bv_1=\Bv
_1(M_k),\ldots
, \Bv_{n-1}=\Bv_{n-1}(M_k)$ and $\Bu_1=\Bu_1(M_k),\ldots, \Bu
_{n-1}=\Bu
_{n-1}(M_k)$ be the orthogonal systems of left and right singular
vectors of $M_k$, and let $\lambda_j=\lambda_j(H_k')= \frac{1}{n}
\sigma
_j^2(M_k)$ be the associated eigenvalues, one has
\[
\Ba_k^\ast\bigl(H_k'-zI
\bigr)^{-1}\Ba_k = \sum_{1\le j\le n-1}
\frac{|\Ba
_k^\ast
\Bv_j|^2}{\lambda_j-z}.
\]

Thus
\[
\Im \bigl(\Ba_k^\ast\bigl(H_k'-zI
\bigr)^{-1}\Ba_k \bigr) \ge\eta\sum
_{1\le
j\le
n-1} \frac{|\Ba_k^\ast\Bv_j|^2}{\eta^2 + |\lambda_j-x|^2}.
\]

We conclude from \eqref{eqnk} that
\[
\sum_{1\le j\le n-1} \frac{|\Ba_k^\ast\Bv_j|^2}{\eta^2 + |\lambda
_j-x|^2}\ll\frac{1}{C\eta}.
\]

Note that $\Ba_k^\ast\Bv_j$ can be written as
\[
\Ba_k^\ast\Bv_j=\frac{\sigma_j(M_k)}{n}
\row_k \Bu_j.
\]

Next, from the Cauchy interlacing law, one can find an interval
$L\subset\{1,\ldots,n-1\}$ of length
\[
|L|\gg C\eta n
\]
such that $\lambda_j\in L$. We conclude that
\[
\sum_{j\in L} \frac{\sigma_j^2}{n^2}|\row_k
\Bu_j|^2 \ll\frac
{\eta}{C}.
\]

Since $\lambda_j\in J$, one has $\sigma_j=\Theta(\sqrt{n})$, and thus
\[
\sum_{j\in L} |\row_k
\Bu_j|^2 \ll\frac{\eta n}{C}.\vadjust{\goodbreak}
\]

The LHS can be written as $\|\pi_V(\row_k^\ast)\|^2$, where $V$ is the
span of the eigenvectors $\Bu_j$ for $j\in L$, and $\pi_V(\cdot)$ is the
projection onto $V$. But from Talagrand's inequality for distance
(Lemma~\ref{lemmadistance} below), we see that this quantity is $\gg
\eta n$ with very high probability, giving the desired contradiction.

\begin{lemma}\label{lemmadistance} Assume that $V\subset\C^n$ is a
subspace of dimension $\dim(V)=d \le n-10$. Let $\Bf$ be a fixed vector
(whose coordinates may depend on $n$). Let $\By=(0,y_2,\ldots,y_n)$,
where $y=\tilde{y}_i-1$ and $\tilde{y}_i$ are i.i.d. copies of
$\tilde
{y}$ defined from~\eqref{eqntilde{y}}. Let $\sigma=\Theta(1)$ denote
the standard deviation of $\tilde{y}$ and $K=10\log n$ denote the upper
bound of $\tilde{y}$, and then for any $t>0$ we have
\[
\P_{\By} \bigl(\pi_V(\By+\Bf)\ge\sqrt{2}\sigma\sqrt
{d}/2-O(K)-t \bigr) \ge1-O \biggl(\exp\biggl(-\frac{t^2}{16K^2}\biggr) \biggr).
\]
\end{lemma}

We now give a proof of Lemma~\ref{lemmadistance}. It is clear that the
function $(y_2,\ldots,\break y_n)\mapsto\pi_V(\By+\Bf)$ is convex and
1-Lipschitz. Thus by Theorem~\ref{lemmaTalagrand} we have
%
%
\begin{equation}
\label{eqnmedium} \P_\By \bigl( \bigl|\pi_V(\By+\Bf)-M\bigl(
\pi_V(\By+\Bf)\bigr)\bigr|\ge t \bigr) =O \bigl(\exp \bigl(-16t^2/K^2
\bigr) \bigr).
\end{equation}

Hence, it is implied that
%
%
\begin{eqnarray}
\label{eqnmedian} &&\P_{\By,\By'} \bigl(\bigl|\pi_V(\By+\Bf)+
\pi_V\bigl(\By'+\Bf\bigr)-2M\bigl(\pi _V(
\By+\Bf )\bigr)\bigr|\le2t \bigr)\nonumber\\
&&\qquad= \bigl(1-O\bigl(\exp\bigl(-16t^2/K^2
\bigr)\bigr) \bigr)^2
\\
&&\qquad=1-O \bigl(\exp\bigl(-16t^2/K^2\bigr) \bigr),\nonumber
\end{eqnarray}
where $\By'$ is an independent copy of $\By$.

On the other hand, by the triangle inequality
\[
\pi_V(\By+\Bf)+\pi_V\bigl(\By'+f\bigr)
\ge\pi_V\bigl(\By-\By'\bigr).
\]

Applying Talagrand's inequality once more for the random vector $\By
-\By
'$ (see, e.g., \cite{TVlocal}, Lemma~68), we see that
\[
\P_{\By,\By'} \bigl(\bigl|\pi_V\bigl(\By-\By'\bigr)-
\sqrt{2}\sigma\sqrt {d}\bigr|\ge t \bigr)=O \bigl(\exp\bigl(-t^2/16K^2
\bigr) \bigr).
\]

Thus,
\[
\P_{\By,\By'} \bigl(\pi_V(\By)+\pi_V\bigl(
\By'\bigr) \ge\sqrt{2}\sigma \sqrt{d} - t \bigr)=1-O \bigl(\exp
\bigl(-t^2/16K^2\bigr) \bigr).
\]

By comparing with \eqref{eqnmedian}, we deduce that
\[
M\bigl(\pi_V(\By+\Bf)\bigr)\ge\sqrt{1/2}\sigma\sqrt{d}-O(K).
\]

Substituting this bound back into \eqref{eqnmedian}, we obtain the
one-sided estimate as desired.
\end{appendix}

\section*{Acknowledgments} The author is grateful to M.~Meckes for
pointing out references \cite{MM} and \cite{P} and to A.~Guionnet for a
helpful e-mail exchange regarding Lemma~\ref
{lemmaconcentrationgeneral}. He is particularly thankful to
R.~Pemantle and V.~Vu for helpful discussions and enthusiastic
encouragement.\vadjust{\goodbreak}


%



\printaddresses

\end{document}